\documentclass[oneside]{article}
\usepackage{amsmath,amssymb,amsthm,graphicx,epsfig,subfigure,float,url}
\usepackage[colorlinks=true]{hyperref}
\usepackage{pdfsync}
\usepackage{color}
\usepackage{esint}
\usepackage{verbatim}
\usepackage[applemac]{inputenc}
\usepackage[usenames,dvipsnames]{pstricks}
\usepackage{pst-grad} % For gradients
\usepackage{pst-plot} % For axes
\usepackage{empheq}
\usepackage{tikz}%Code Matlab
\usepackage{mathrsfs}

\usepackage{listings}%Code Matlab
\usepackage[usenames,dvipsnames]{xcolor}
\hypersetup{ citecolor=gray, linkcolor=purple, urlcolor=blue}

\def\ls{{\lesssim}}

\def\R{\textrm{I\kern-0.21emR}}
\def\N{\textrm{I\kern-0.21emN}}

\def\O{{\Omega}}
\def\n{{\nabla}}
\def\tmm{{\theta_{m,\mu}}}
\def\pmm{{p_{m,\mu}}}
\def\p{{\varphi}}
\def\phimm{{\varphi_{m,\mu}}}
\def\etam{{\eta_{1,m}}}
\def\etakm{{\eta_{k,m}}}
\def\etakkm{{\eta_{k+1,m}}}
\def\hetakkm{{\hat\eta_{k+1,m}}}

\def\betakkm{{\beta_{k+1,m}}}

\def\detakm{{\dot{\eta}_{k,m}}}

\def\ekm{{\eta_{k,m}}}

\def\lp{{L^p(\O)}}

\def\hetam{{\hat{\eta}_{1,m}}}
\def\hetakm{{\hat{\eta}_{k,m}}}

\def\dhetakkm{{\dot{\hat\eta}_{k+1,m}}}

\def\hetakm{{\hat{\eta}_{k,m}}}
\def\betakm{{\beta_{k,m}}}
\def\hetam{{\hat{\eta}_{1,m}}}
\def\dhetam{{\dot{\hat{\eta}}_{1,m}}}

\def\betam{{\beta_{1,m}}}

\def\e{{\varepsilon}}

\def\D{{\Delta}}
\def\l{{L^2(\O)}}
\def\li{{L^\infty(\O)}}
\def\dtm{{\dot{\theta}_{m,\mu}}}

\renewcommand{\geq}{\geqslant}
\renewcommand{\leq}{\leqslant}
\renewcommand{\red}[1]{\textcolor{black}{#1}}

\newtheorem{theorem}{Theorem}
\newtheorem*{thmnonumbering}{Theorem}

\newtheorem{proposition}{Proposition}
\newtheorem{corollary}{Corollary}
\newtheorem{definition}{Definition}
\newtheorem{lemma}{Lemma}
\newtheorem{Remark}{Remark}

\theoremstyle{definition}
\theoremstyle{definition}\newtheorem{remark}{Remark}

\title{Optimal location of resources maximizing the total population size in logistic models\footnote{The authors were partially supported by the Project ``Analysis and simulation of optimal shapes - application to lifesciences'' of the Paris City Hall. I. Mazari and Y. Privat were partially supported by the ANR Project ``SHAPe Optimization - SHAPO''.}}

\author{Idriss Mazari\footnote{Sorbonne Universit\'es, UPMC Univ Paris 06, UMR 7598, Laboratoire Jacques-Louis Lions, F-75005, Paris, France (\texttt{idriss.mazari@sorbonne-universite.fr}).}
       \and Gr\'egoire Nadin\footnote{ CNRS, Sorbonne Universit\'es, UPMC Univ Paris 06, UMR 7598, Laboratoire Jacques-Louis Lions, F-75005, Paris, France (\texttt{gregoire.nadin@sorbonne-universite.fr})}
	\and Yannick Privat\footnote{IRMA, Universit\'e de Strasbourg, CNRS UMR 7501, 7 rue Ren\'e Descartes, 67084 Strasbourg, France ({\tt yannick.privat@unistra.fr}).}
	}

\date{}

\begin{document}

\maketitle

\begin{abstract}
In this article, we consider a species whose population density solves the steady diffusive logistic equation in a heterogeneous environment modeled with the help of a spatially non constant coefficient standing for a resources distribution. We address the issue of maximizing the total population size with respect to the resources distribution, considering some uniform pointwise bounds as well as prescribing the total amount of resources. By assuming the diffusion rate of the species large enough, we prove that any optimal configuration is bang-bang (in other words an extreme  point of the admissible set) meaning that this problem can be recast as a shape optimization problem, the unknown domain standing for the resources location. In the one-dimensional case, this problem is deeply analyzed, and for large diffusion rates, all optimal configurations are exhibited. This study is completed by several numerical simulations in the one dimensional case.
\end{abstract}

\medskip

%\tableofcontents

%\noindent{\bf Acknowledgment.}
\tableofcontents

%%%%%%%%%%%%%%%%%%%%%%%%%%%%%%%%%%%%%%%%%%%Introduction
\section{Introduction}

%%Notations et historiques
\subsection{Motivations and state of the art}\label{sec:motiv}
In this article, we  investigate an optimal control problem arising in population dynamics. Let us consider the population density $\tmm$ of a given species evolving in a bounded and connected domain $\O$ in $\R^n$ with $n\in \N^*$, having a $\mathscr C^2$ boundary. In what follows, we will assume that $\tmm$ is the positive solution of the  steady logistic-diffusive equation (denoted  \eqref{LDE} in the sequel)  which writes 
\begin{equation}\tag{LDE}\label{LDE} \left\{
\begin{array}{ll}
\mu \D \tmm(x)+(m(x)-\tmm(x))\tmm(x) =0 & x\in \O,\\
\frac{\partial \tmm}{\partial \nu} =0 &x\in \partial \O,
%\\\tmm&\geq &0,
\end{array}
\right.\end{equation} 
where $m\in L^\infty(\O)$ stands for the resources distribution and $\mu>0$ stands for the {dispersal ability}  of the species, also called {\it diffusion rate}. From a biological point of view, the real number $m(x)$ is the local intrinsic growth rate of species at location $x$ of the habitat $\O$ and can be seen as a measure of the resources available at $x$.

As will be explained below, we will only consider non-negative resource distributions.$m$, i.e such that $m\in L^\infty_+(\O)=\left\{m\in L^\infty(\O)\, , m\geq 0 \text{ a.e}\right\}$. 
In view of investigating the influence of spatial heterogeneity on the model, we consider the optimal control problem  of maximizing the functional  
$$
{\mathcal F_\mu}:L_+^\infty(\O)\ni m\mapsto \fint_\O \tmm,
$$
{where the notation $\fint$ denotes the average operator, in other words $\fint_\O f=\frac1{|\O|}\int_\O f$.}
The functional $F$ stands for the total population size, in order to further our understanding of spatial heterogeneity on population dynamics. 
%Here, for computational and notational simplicity, we  use the average  of a function, i.e $\fint_\O f=\frac1{|\O|}\int_\O f$ rather than the integral $\int_\O f$. 
\\In the framework of population dynamics, the density $\theta_{m,\mu}$ solving Equation \eqref{LDE} can be interpreted as a steady state associated to the following evolution equation
\begin{equation}\tag{LDEE}\label{LDEE} \left\{
\begin{array}{ll}
 \frac{\partial u}{\partial t}(t,x) = \mu \D u(t,x)+u(t,x)(m(x)-u(t,x)) & t>0, \ x\in \O\\
\frac{\partial u}{\partial \nu}(t,x)=0& t>0, \ x\in \partial \O\\ 
u(0,x)=u^0(x)\geq 0\, , u^0\neq 0 & x\in \O \\
%u&\geq &0.
\end{array}
\right.\end{equation}
modeling the spatiotemporal behavior of a population density $u$ in a domain $\Omega$ with the spatially heterogeneous resource term $m$.

The pioneering works by Fisher \cite{Fisher}, Kolmogorov-Petrovski-Piskounov \cite{KPP} and Skellam \cite{SkellamRandom} on the logistic diffusive equation were mainly concerned with the spatially homogeneous case. Thereafter, many authors investigated the influence of spatial heterogeneity on population dynamics and species survival. In \cite{Shigesada3},  propagation properties in a patch model environment are studied. In \cite{Shigesada2}, a spectral condition for species survival in heterogeneous environments has been derived, while \cite{Shigesada1} deals with the influence of fragmentation and concentration of resources on population dynamics. These works were followed by \cite{BHR} dedicated to an optimal design problem, that will be commented in the sequel.

Investigating existence and uniqueness properties of solutions for the two previous equations as well as their regularity properties boils down to the study of spectral properties for the linearized operator 
$$
\mathcal L:\mathcal{D}(\mathcal{L})\ni f\mapsto \mu\Delta f+mf,
$$
where  the domain of $\mathcal L$ is $\mathcal{D}(\mathcal{L})=\{f\in L^2(\O)\mid \Delta f\in L^2(\O)\}$ and of its first eigenvalue $\lambda_1(m,\mu)$,  characterized by the Courant-Fischer formula
\begin{equation}\label{CFformulaL1}
\lambda_1(m,\mu):=\underset{f\in W^{1,2}(\Omega),\int_\Omega f^2=1}\sup \left\{-\mu \int_\Omega |\nabla f|^2+\int_\Omega mf^2\right\}.
\end{equation}
Indeed, the positiveness of $\lambda(m,\mu)$ is a sufficient condition ensuring the well-posedness of equations \eqref{LDEE} and \eqref{LDE} (\cite{BHR}). Then, Equation \eqref{LDE} has a unique positive solution $\tmm\in W^{1,2}(\O)$. Furthermore, for any $p\geq 1$, $\tmm$ belongs to $W^{2,p}(\O)$, and there holds
\begin{equation}\label{metz:1944}
0<\underset{\overline \O}\inf\, \tmm\leq \tmm\leq \Vert m\Vert _{L^\infty(\O)}.
\end{equation}
Moreover, the steady state $\tmm$ is globally asymptotically stable: for any $u_0\in W^{1,2}(\O)$ such that $u_0\geq0$ a.e. in $\O$ and $u_0\neq 0$ , one has
$$
\Vert u(t,\cdot)-\tmm\Vert_{L^\infty(\O)}\underset{t\rightarrow +\infty}\longrightarrow 0.
$$
where $u$ denotes the unique solution of \eqref{LDEE} with initial state $u_0$ (belonging to $L^2(0,T;W^{1,2}(\O))$ for every $T>0$).

The importance of $\lambda_1(m,\mu)$ for stability issues related to population dynamics models was first noted in simple cases by Ludwig, Aronson and Weinberger \cite{AronsonLudwigWeinbergerSpatial}. Let us mention \cite{DockeryHutsonMischaikowPernarowskiEvolution} where the case of diffusive Lotka-Volterra equations is investigated. 

To guarantee that $\lambda_1(m,\mu)>0$, it is enough to work with distributions of resources $m$ satisfying the assumption
\begin{equation}\label{H1}\tag{H1}
m\in L^\infty_+(\O)\quad \text{where}\quad L^\infty_+(\O)=\left\{m\in L^\infty(\O), \ \ \int_\O m>0\right\}.
\end{equation}
Note that the issue of maximizing this principal eigenvalue was addressed for instance in \cite{JhaPorru,KaoLouYanagida,LamboleyLaurainNadinPrivatProperties,MR2281509,Roques-Hamel}.

In the survey article \cite{LouSome}, Lou suggests the following problem: the parameter $\mu>0$ being fixed, which weight $m$ maximizes the total population size among all uniformly bounded elements of $L^\infty(\O)$?

In this article, we aim at providing partial answers to this issue, and more generally new results about  the influence of the spatial heterogeneity  $m(\cdot)$ on the total population size.

For that purpose, let us introduce the total population size functional, defined for a given $\mu>0$ by 
\begin{equation}\label{def:Fmu}
{\mathcal F_\mu}: L^\infty_+(\O)\ni m\longmapsto \fint_\O \tmm,
\end{equation}
where $ \tmm$ denotes the solution of equation \eqref{LDE}. 

\medskip 

Let us mention several previous works dealing with the maximization of the total population size functional. It is shown in \cite{LouMigration} that, among all weights $m$ such that $\fint_\O m=m_0$, there holds
  ${\mathcal F_\mu}(m)\geq {\mathcal F_\mu}(m_0)=m_0$; this inequality is strict  whenever $m$ is nonconstant. 
Moreover, it is also shown that the problem of maximizing ${\mathcal F_\mu}$ over $ L^\infty_+(\O)$ has no solution.
  
   \begin{remark}
 \label{rk:min}The fact that $m\equiv m_0$ is a minimum for ${\mathcal F_\mu}$ among the resources distributions $m$ satisfying $\fint_\O m=m_0$  relies on the following observation:
 multiplying \eqref{LDE} by $\frac1\tmm$ and integrating by parts yields
\begin{equation}\label{metz2135}
\mu\fint_\O \frac{|\n \tmm|^2}{\tmm^2}+\fint_\O (m-\tmm)=0.
\end{equation}
and therefore, ${\mathcal F_\mu}(m)=m_0+\mu\int_\O \frac{|\n \tmm|^2}{\tmm^2}\geq m_0={\mathcal F_\mu}(m_0)$  for all $m\in L^\infty_+(\O)$ such that $\fint_\O m=m_0$.   \color{black} It follows that the constant function equal to $m_0$ is a global minimizer of ${\mathcal F_\mu}$ over $ \{m\in L^\infty_+(\O)\, , \fint_\O m=m_0\}.$\color{black}
   \end{remark}

 In the recent article \cite{BaiHeLiOptimization}, it is shown that, when $\O=(0,\ell)$, one has
$$
\forall \mu>0, \ \forall m\in L_+^\infty(\O)\mid m\geq 0\text{ a.e.}, \qquad \fint_\O \tmm\leq 3\fint_\O m.
$$
This inequality is sharp, although the right-hand side is never reached, and the authors exhibit a sequence $(m_k,\mu_k)_{k\in \N}$ such that $\fint_\O \theta_{m_k,\mu_k}/\fint_\O m_k \rightarrow 3$ as $k\to +\infty$, but for such a sequence there holds $\Vert m_k\Vert_{L^\infty(\O)}\rightarrow +\infty$ and $\mu_k\rightarrow 0$ as $k\to +\infty$.
\\In \cite{LouSome}, it is proved that, without $L^1$ or $L^\infty$ bounds on the weight function $m$, the maximization problem is ill-posed. It is thus natural to introduce two parameters $\kappa,m_0>0$, and to restrict our investigation to the class 
\begin{equation}\label{Admissibles}\mathcal M_{m_0,\kappa}(\O):=\left\{m\in L^\infty(\O)\, , 0\leq m\leq \kappa \text{ a.e }\, , \fint_\O m=m_0\right\}.\end{equation}
It is notable that in \cite{DING2010688}, the more general functional $J_B$ defined by 
$$J_B(m)=\int_\O (\theta_{m,\mu}-Bm^2)\qquad  \hbox{ for } B\geq 0$$
is introduced. In the case $B=0$, the authors apply the so-called Pontryagin principle, show the G\^ateaux-differentiability of $J_B$ and carry out numerical simulations backing up the conjecture that maximizers  of $J_0$ over $\mathcal M_{m_0,\kappa}(\O)$ are of bang-bang type.

However, proving this {\it bang-bang} property is a challenge. The analysis of optimal conditions is quite complex, because the sensitivity of the functional "total population size" with respect to $m(\cdot)$ is directly related to the solution of an adjoint state, solving a linearized version of \eqref{LDE}. Deriving and exploiting the properties of optimal configurations therefore requires a thorough understanding of the $\theta_{m,\mu}$ behavior as well as the associated state.  To do this, we are introducing a new asymptotic method to exploit optimal conditions.

We will investigate two properties of the maximizers of the total population size function ${\mathcal F_\mu}$.
\begin{enumerate}
\item \textbf{Pointwise constraints.}
The main issue that will be addressed in what follows is the {\it bang-bang} character of optimal weights $m^*(\cdot)$, in other words, whether $m^*$ is equal to 0 or $\kappa$ almost everywhere. Noting that $\mathcal M_{m_0,\kappa}(\O)$ is a convex set and that bang-bang functions are the extreme points of this convex set, this question rewrites:
\begin{center} \textit{Are the maximizers $m^*$ extreme points of the set $\mathcal{M}_{m_0,\kappa}(\O)$?}\end{center} 
In our main result (Theorem \ref{TheoremePrincipal}) we provide a \textbf{positive answer for large diffusivities}. It is notable that our proof rests upon a well-adapted expansion of the solution $\tmm$ of \eqref{LDE} with respect to the diffusivity $\mu$.

This approach could be considered unusual, since such results are usually obtained by an analysis of the properties of the adjoint state (or switching function). However, since the switching function very implicitly depends on the design variable $m(\cdot)$, we did not obtain this result in this way.
\item \textbf{Concentration-fragmentation.} It is well known that resource concentration (which means that the distribution of resources $m$ decreases in all directions, see Definition \ref{DR} for a specific statement) promotes the survival of species \cite{BHR}. 
On the contrary, we will say that a resource distribution $m=\kappa \chi_E$, where $E$ is a subset of $\O$, is \textit{fragmented} when the $E$ set is disconnected. In the figure \ref{figFC}, $\O$ is a square, and the intuitive notion of \textit{concentration-fragmentation} of resource distribution is illustrated.

\begin{figure}[h]
\begin{center}
\includegraphics[width=5.2cm]{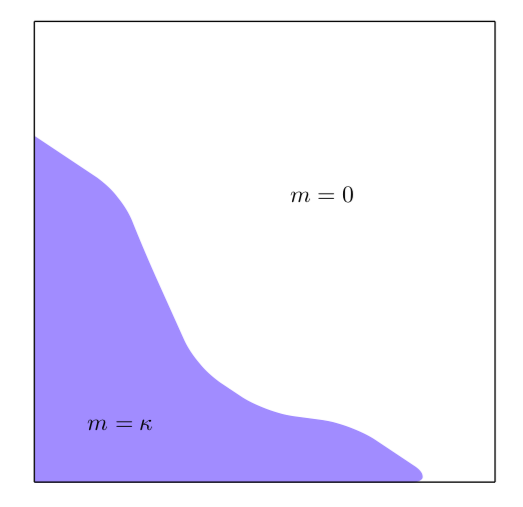}\hspace{1cm}
\includegraphics[width=5.1cm]{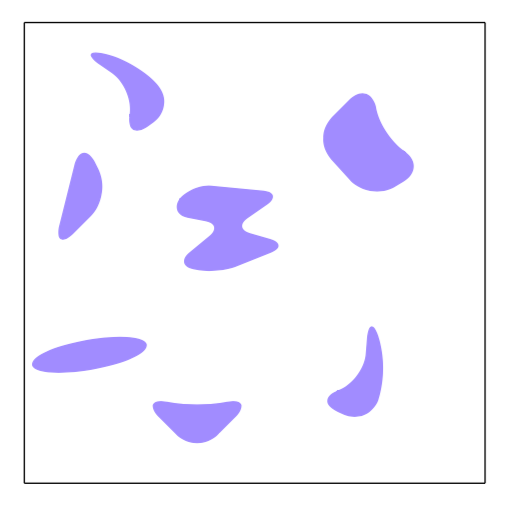}
\caption{$\O=(0,1)^2$. The left distribution 
is  "concentrated" (connected) whereas the right one is fragmented (disconnected). \label{figFC}}
\end{center}
\end{figure}

A natural issue related to qualitative properties of maximizers is thus
\begin{center}
\textit{Are maximizers $m^*$ concentrated? Fragmentated? }
\end{center}
In Theorem \ref{Concentration}, we consider the case of a orthotope shape habitat, and we show that \textbf{concentration occurs for large diffusivities:} if $m_\mu$  maximizes ${\mathcal F_\mu}$ over $\mathcal M_{m_0,\kappa}(\O)$,  the sequence $\{m_\mu\}_{\mu>0}$   strongly converges in $L^1(\O)$ to a concentrated distribution as $\mu \to \infty$,.
\\ In the one-dimensional case, we also prove that {if the diffusivity is large enough, there are only two maximizers}, that are plotted on Fig. \ref{figCDG} (see Theorem \ref{Creneau}).

\begin{figure}[h]
\begin{center}
\includegraphics[width=5cm]{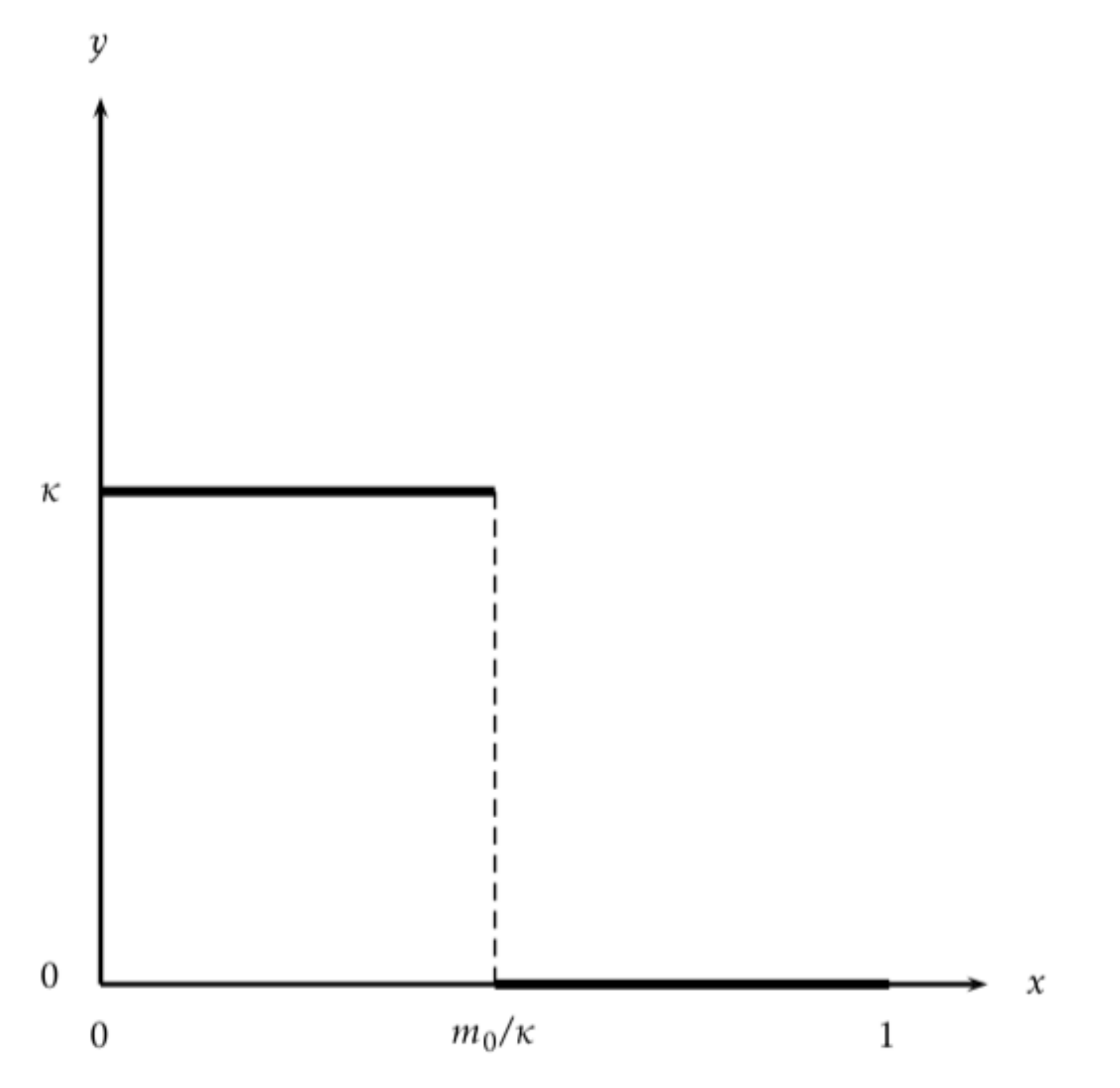}\hspace{1cm}
\includegraphics[width=5.1cm]{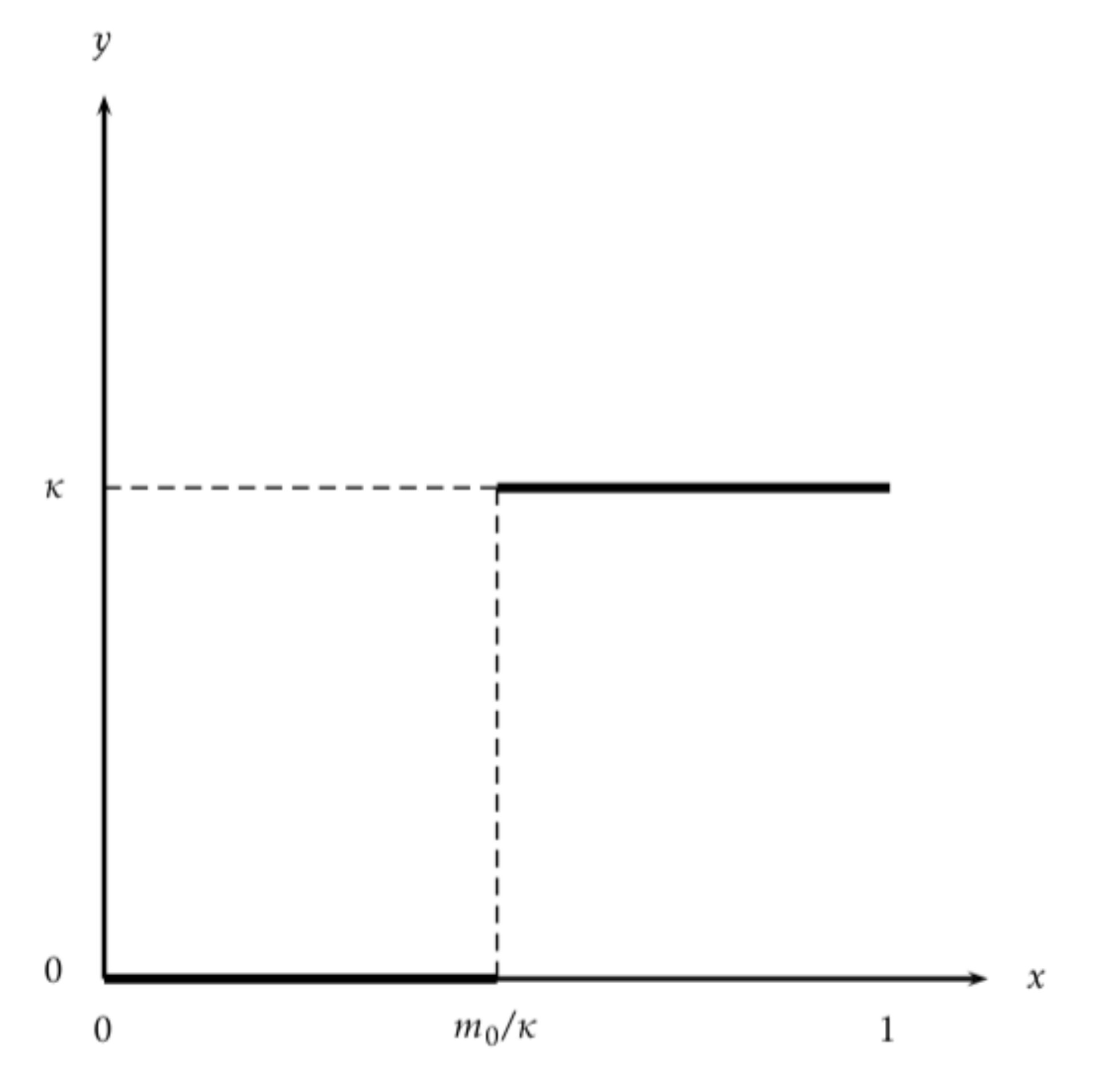}
\caption{$\O=(0,1)$. Plot of the only two maximizers of ${\mathcal F_\mu}$ over $\mathcal{M}_{m_0,\kappa}(\O)$. \label{figCDG}}
\end{center}
\end{figure}
 Finally, {in the one-dimensional case, we obtain a surprising result:} \textbf{fragmentation may be better than concentration for small diffusivities} (see Theorem \ref{fragmentation} and Fig. \ref{figDCCD} below). 

\begin{figure}[h]
\begin{center}
\includegraphics[width=5cm]{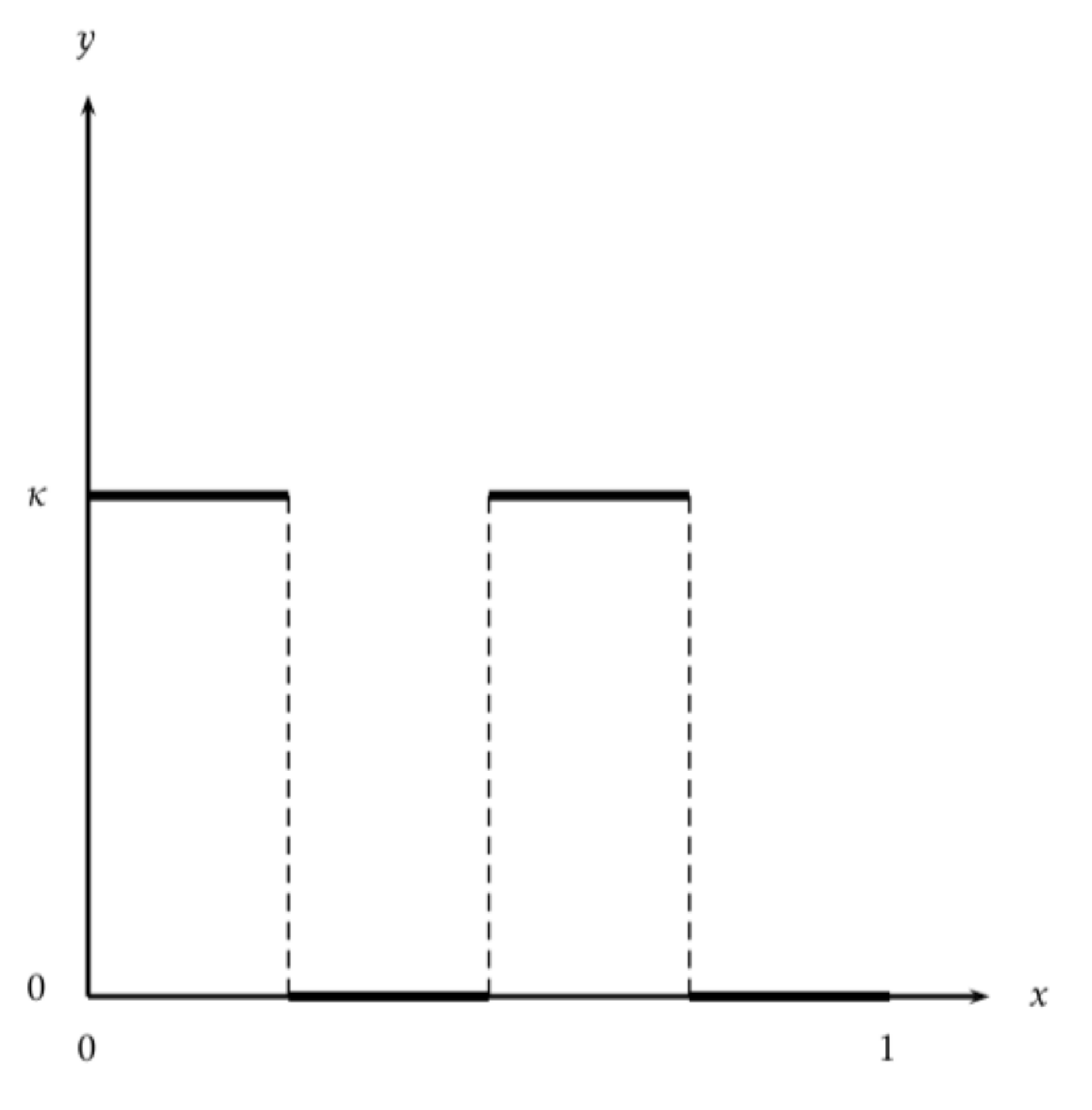}\hspace{1cm}
\includegraphics[width=5.1cm]{creneaudroit.pdf}
\caption{$\O=(0,1)$. A double crenel (on the left) is better than a single one (on the right). \label{figDCCD}}
\end{center}
\end{figure}

 This is surprising because in many problems of optimizing the logistic-diffusive equation, it is expected that the best disposition of resources will be concentrated.
 
\end{enumerate}

\subsection{Main results}

In the whole article, the notation $\chi_I$ will be used to denote the characteristic function of a measurable subset $I$ of $\R^n$, in other words, the function equal to 1 in $I$ and 0 elsewhere.

For the reasons mentioned in Section \ref{sec:motiv}, it is biologically relevant to consider the class of admissible weights
$
\mathcal{M}_{m_0,\kappa}(\O)$ defined by \eqref{Admissibles},
where $\kappa>0$ and $m_0>0$ denote two positive  parameters such that $m_0<\kappa$ (so that this set is nontrivial).

We will henceforth consider the following optimal design problem.

\begin{center}
\fbox{
\begin{minipage}{15cm}
\begin{center}
\begin{quote}
{\bf Optimal design problem. }\textit{
Fix $n\in \N^*$, $\mu>0$, $\kappa>0$, $m_0\in (0,\kappa)$ and let $\O$ be a bounded connected domain of $\R^n$ having a $\mathscr C^2$ boundary.
We consider the optimization problem
\begin{equation}\tag{$\mathcal{P}_\mu^n$}\label{TSOM}
\sup_{m\in \mathcal{M}_{	\kappa,m_0}(\O)} {\mathcal F_\mu}(m).
\end{equation} }
\end{quote}
\end{center}\end{minipage}}
\end{center}

As  will be highlighted in the sequel, the existence of a maximizer follows from a direct argument.
We will thus investigate the  qualitative properties of maximizers  described in the previous section (bang-bang character, concentration/fragmentation phenomena).
%Quid of the uniqueness? 

For the sake of readability, almost all the proofs are postponed to Section \ref{sec:proofMainres}.

\medskip

Let us stress that the bang-bang character of maximizer is of practical interest in view of spreading resources in an optimal way. Indeed, in the case where a maximizer $m^*$ writes $m^*=\kappa \chi_E$, the total size of population is maximized by locating all the resources on $E$. 
\subsubsection{First property}
We start with a preliminary result related to the saturation of pointwise constraints for Problem \eqref{TSOM}, valid for all diffusivities $\mu$.  It is obtained by exploiting the first order optimality conditions for Problem \eqref{TSOM}, written in terms of an adjoint state. 

\begin{proposition}\label{SaturationContraintePonctuelle}
Let $n\in \N^*$, $\mu>0$, $\kappa>0$, $m_0\in (0,\kappa)$. Let $m^*$ be a solution of Problem \eqref{TSOM}. Then, the set $\{m=\kappa\}\cup \{m=0\}$ has a positive measure
\end{proposition}

\subsubsection{The bang-bang property holds for large diffusivities}

For large values of $\mu$, we will prove that the variational problem can be recast in terms of a shape optimization problem, as underlined in the next results.

\begin{theorem}\label{TheoremePrincipal}
Let $n\in \N^*$, $\kappa>0$, $m_0\in (0,\kappa)$. There exists a positive number $\mu^*=\mu^*(\O,\kappa,m_0)$ such that, for every $\mu\geq \mu^*$, the functional ${\mathcal F_\mu}$ is strictly convex. As a consequence, for $\mu\geq \mu^*$, {any} maximizer of ${\mathcal F_\mu}$ over $\mathcal{M}_{m_0,\kappa}(\O)$ (or similarly any solution of Problem \eqref{TSOM}) is moreover of {\it bang-bang} type\footnote{In other words, it is,an element of $\mathcal{M}_{m_0,\kappa}(\O)$ equal a.e. to 0 or $\kappa$ in $\O$.}.
\end{theorem}

We emphasize that the proof of Theorem \ref{TheoremePrincipal}  is quite original, since it does not rest upon the exploitation of adjoint state properties, but upon the use of a power series expansions in the diffusivity $\mu$ of the solution $\tmm$ of \eqref{LDE}, as well as their derivative with respect to the design variable m. In particular, this expansion is used to prove that, if $\mu$ is large enough, then the function ${\mathcal F_\mu}$ is strictly convex. Since the extreme points of $\mathcal M_{m_0,\kappa}(\O)$ are bang-bang resources functions, the conclusion readily follows.
\\This theorem justifies that Problem \eqref{TSOM} can be recast as a shape optimization problem. Indeed, every maximizer $m^*$ is of the form $m^*=\kappa\chi_{E}$ where $E$ is a measurable subset such that $|E|=m_0|\O|/\kappa$. 
\begin{remark} {
We can rewrite this result in terms of shape optimization, by considering as main unknown the subset $E$ of $\O$ where resources are located: indeed, under the assumptions of Theorem \ref{TheoremePrincipal}, there exists a positive number $\mu^*=\mu^*(\O,\kappa,m_0)$ such that, for every $\mu\geq \mu^*$, the shape optimization problem
\begin{equation}\label{SOP}
\sup_{E\subset \O, \ |E|=m_0|\O|/\kappa} {\mathcal F_\mu}(\kappa\chi_{E}),
\end{equation}
where the supremum is taken over all measurable subset $E\subset \O$ such that  $|E|=m_0|\O|/\kappa$, has a solution.
We underline the fact that, for such a shape functional, which is ``non-energetic'' in the sense that the solution of the PDE involved cannot be seen as a minimizer of the same functional, proving the existence of maximizers is usually intricate.}
\end{remark}
\subsubsection{Concentration occurs for large diffusivities}
In this section, we state two results suggesting concentration properties for the solutions of Problem \eqref{TSOM} may hold for large diffusivities.
\\For that purpose, let us introduce the function space
\begin{equation}\label{EspaceNul}
X:=W^{1,2}(\O)\cap \left\{u\in W^{1,2}(\O)\, , \fint_\O u=0\right\}
\end{equation}
and the energy functional 
\begin{equation}\label{def:Em}
\mathcal E_m:X\ni u\mapsto \frac12  \fint_\O |\n u|^2-m_0\fint_\O mu.
\end{equation}
\begin{theorem}\label{Concentration}[$\Gamma$-convergence property]
\begin{enumerate}
\item Let $\O$ be a domain with a  $\mathscr C^2$ boundary. For any $\mu>0$, let  $m_\mu$ be a solution of Problem \eqref{TSOM}. Any $L^1$ closure point of $\{m_\mu\}_{\mu>0}$ as $\mu \to \infty$ is a solution of the optimal design problem 
\begin{equation}\label{Pb:Lim}\min_{m\in \mathcal M_{m_0,\kappa}(\O)}\min_{u\in X} \mathcal E_m(u).\end{equation}
\item In the case of a two dimensional orthotope $\O=(0;a_1)\times (0;a_2)$, any solution of the asymptotic optimization problem \eqref{Pb:Lim} decreases in every direction.
\end{enumerate}
\end{theorem}
 As will be clear in the proof, this theorem is a $\Gamma$-convergence property.

In the one-dimensional case, one can refine this result by showing that, for $\mu$ large enough, the maximizer is a step function. 
\begin{theorem}\label{Creneau}
Let us assume that $n=1$ and $\O=(0,1)$. Let  $\kappa>0$, $m_0\in (0,\kappa)$. There exists $\hat \mu>0$ such that, for any $\mu\geq \hat \mu$, any solution $m$ of Problem \eqref{TSOM} is equal a.e. to either $\tilde m$ or $\tilde m(1-\cdot)$, where $\tilde m=\kappa \chi_{(1-\ell,1)}$ and $\ell=m_0/\kappa$. 
\end{theorem}
\subsubsection{Fragmentation may occur for small diffusivities}
Let us conclude by underlining that the statement of Theorem \ref{Creneau} cannot be true for all $\mu>0$. Indeed, we provide an example in Section \ref{sec:1D} where a double-crenel growth rate gives a larger total population size than the simple crenel $\tilde m$ of Theorem \ref{Creneau}.

\begin{theorem}\label{fragmentation}The function $\tilde m=\kappa \chi_{(1-\ell,1)}$ (and  $\tilde m(1-\cdot)=\kappa\chi_{(0,\ell)}$) does not solve Problem \eqref{TSOM} for small values of $\mu$.
More precisely, if we extend $\tilde m$ outside of $(0,1)$ by periodicity, there exists $\mu>0$ such that
$$
 {\mathcal F_\mu}\big(\tilde m(2~\cdot )\big) > {\mathcal F_\mu}(\tilde m).
$$
  \end{theorem}
This result  is quite unusual. For the optimization of the first eigenvalue $\lambda_1(m,\mu)$ defined by \eqref{CFformulaL1} with respect to $m$ on the interval $(0;1)$
$$\sup_{m\in \mathcal M((0;1))} \lambda_1(m,\mu)$$ we know (see \cite{BHR}) that the only solutions are $\tilde m$ and $\tilde m(1-\cdot)$,  for any $\mu$.

It is notable that the following result is a byproduct of Theorem \ref{fragmentation} above.
 \begin{corollary}
 There exists $\mu>0$ such that the problems
 $$\sup_{m\in \mathcal M((0;1))} \lambda_1(m,\mu)$$ and 
 $$\sup_{m\in \mathcal M((0;1))} {\mathcal F_\mu}(m)$$ do not have the same maximizers.
 \end{corollary}
For further comments on the relationship between the main eigenvalue and the total size of the population, we refer to \cite{Mazari}, where an asymptotic analysis of the main eigenvalue (relative to $\mu$ as $\mu \to +\infty$) is performed, and the references therein.

We conclude this section by mentioning the recent work \cite{NagaharaYanagida}, that was reported to us when we wrote this article. They show that, if we assume the optimal distribution of regular resources (more precisely, Riemann integrable), then it is necessarily of {\it bang-bang} type. Their proof is based on a perturbation argument valid for all $\mu>$0. However, proving such regularity is generally quite difficult. Our proof, although it is not valid for all $\mu$, is not based on such a regularity assumption, but these two combined results seem to suggest that all maximizers of this problem are of {\it bang-bang} type.
\color{black}

\subsection{Tools and notations}\label{sec:toolsNota}
In this section, we gather some useful tools we will use to prove the main results.

\paragraph{Rearrangements of functions and principal eigenvalue.}
Let us first recall several monotonicity and regularity related to the principal eigenvalue of the operator $\mathcal L$.
\begin{proposition}\cite{DockeryHutsonMischaikowPernarowskiEvolution}
Let $m\in L^\infty_+(\O)$ and $\mu>0$.
\begin{itemize}
\item[(i)] The mapping $\R_+^*\ni \mu\mapsto \lambda_1(m,\mu)$ is continuous and non-increasing.
\item[(ii)] If $m\leq m_1$, then $\lambda_1(m,\mu)\leq \lambda_1(m_1,\mu)$, and the equality is true if, and only if $m=m_1$ a.e. in $\O$.
\end{itemize}
\end{proposition}

In the proof of Theorem \ref{Creneau}, we will  use rearrangement inequalities at length. Let us briefly recall the notion of decreasing rearrangement. 
\begin{definition}
For a given function $b\in L^2(0,1)$, one defines its monotone decreasing  (resp. monotone increasing) rearrangement $b_{dr}$ (resp. $b_{br}$) on $(0,1)$ by $b_{dr}(x)=\sup \{c\in \R \mid x\in \Omega_c^*\}$, where $\Omega_c^*=(1-|\Omega_c|,1)$ with $\Omega_c=\{b>c\}$ (resp. $b_{br}(\cdot)=b_{dr}(1-\cdot)$).
\end{definition}
The functions $b_{dr}$ and $b_{br}$ enjoy nice properties. In particular, the Poly\`a-Szeg{o} and Hardy-Littlewood inequalities allow to compare integral quantities depending on $b$, $b_{dr}$, $b_{br}$ and their derivative. 

\begin{thmnonumbering}[\cite{MR810619,LiebLoss}]
Let $u$ be a non-negative and measurable function.
\begin{itemize}
\item[(i)] If $\psi$ is any measurable function from $\R_+$ to $\R$, then 
$$
\int_0^1 \psi(u)=\int_0^1 \psi(u_{dr})=\int_0^1 \psi(u_{br})\quad \text{(equimeasurability)};
$$
\item[(ii)] If $u$ belongs to $W^{1,p}(0,1)$ with $1\leq p$, then 
$$
\int_0^1(u')^p\geq \int_0^1(u_{br}')^p=\int_0^1(u_{dr}')^p \quad \text{(P\'olya inequality)};
$$
\item[(iii)] If $u$, $v$ belong to $L^2(0,1)$, then
$$
\int_0^1 uv\leq \int_0^1 u_{br}v_{br}=\int_0^1 u_{dr}v_{dr}  \quad \text{(Hardy-Littlewood inequality)};
$$
\end{itemize}
\end{thmnonumbering}
The equality case in the Poly\`a-Szego inequality is the object of the Brothers-Ziemer theorem (see e.g. \cite{FeroneVolpicelli}).

\paragraph{Symmetric decreasing functions}\label{Ptes}
In higher dimensions, in order to highlight concentration phenomena, we will use another notion of symmetry, namely monotone symmetric rearrangements that are extensions of monotone rearrangements in one dimension. Here, $\O$ denotes the $n$-dimensional orthotope $\prod_{i=1}^n(a_i,b_i)$. 
\begin{definition}\label{DR}
For a given function $b\in L^1(\O)$, one defines its symmetric decreasing  rearrangement $b_{sd}$  on $\O$ as follows: first fix the $n-1$ variables $x_2,\dots,x_n$. Define $b_{1,sd}$ as the monotone decreasing rearrangement of $x\mapsto b(x,x_2,\dots,x_n)$. Then fix $x_1,x_3,\dots,x_n$ and define $b_{2,sd}$ as the monotone decreasing rearrangement of $x\mapsto b_{1,sd}(x_1,x,\dots,x_n)$. Perform such monotone decreasing rearrangements successively. The resulting function is the symmetric decreasing rearrangement of $b$.
\\We define the symmetric increasing rearrangement in a direction $i$ a similar fashion and write it $b_{i,id}$.  Note that, in higher dimensions, the definition of decreasing rearrangement strongly depends on the order in which the variables are taken.\color{black}
\end{definition}

Similarly to the one-dimensional case, the P\'olya-Szego and Hardy-Littlewood inequalities allow us to compare integral quantities.

\begin{thmnonumbering}[\cite{BHR,BLR}]\label{RSD}
Let $u$ be a non-negative and measurable function defined on a box $\O=\prod_{i=1}^n(0;a_i)$.
\begin{itemize}
\item[(i)] If $\psi$ is any measurable function from $\R_+$ to $\R$, then 
$$
\int_\O \psi(u)=\int_\O\psi(u_{sd})=\int_\O \psi(u_{sd})\quad \text{(equimeasurability)};
$$
\item[(ii)] If $u$ belongs to $W^{1,p}(\O)$ with $1\leq p$, then, for every $i\in \N_N$, $[a_i;b_i]=\omega_{1,i}\cup \omega_{2,i}\cup \omega_{3,i}$, where the map $(x_1,\dots,x_{i-1},x_{i+1},\dots,x_N)\mapsto u(x_1,\dots,x_n)$ is decreasing if $x_i\in \omega_{1,i}$, increasing if $ x_i\in \omega_{2,i}$ and constant if $x_i\in \omega_{3,i}$.
$$
\int_\O|\n u|^p\geq \int_\O|\n u_{sd}|^p \quad \text{(P\'olya inequality)};
$$
Furthermore, if, for any $i\in \N_n$, $\int_\O|\n u|^p= \int_\O|\n u_{i,sd}|^p $ then there exist three measurable subets $\omega_{i,1}\, , \omega_{i,2}$ and $\omega_{i,3}$ of $(0;a_i)$ such that
\begin{enumerate}
\item $(0;a_i)=\omega_{i,1}\cup  \omega_{i,2}\omega_{i,3}$,
\item $u=u_{i,bd}$ on $\prod_{k=1}^{i-1}(0;a_k)\times \omega_{i,1}\times \prod_{k=i+1}^n (0;a_k)$,
\item $u=u_{i,id}$ on $\prod_{k=1}^{i-1}(0;a_k)\times \omega_{i,2}\times \prod_{k=i+1}^n (0;a_k)$,
\item $u=u_{i,bd}=u_{i,id}$ on $\prod_{k=1}^{i-1}(0;a_k)\times \omega_{i,3}\times \prod_{k=i+1}^n (0;a_k)$.
\end{enumerate}
\item[(iii)] If $u$, $v$ belong to $L^2(\O)$, then
$$
\int_\O uv\leq \int_\O u_{sd}v_{sd}  \quad \text{(Hardy-Littlewood inequality)};
$$
\end{itemize}
\end{thmnonumbering}

\paragraph{Poincar\'e constants and elliptic regularity results.} 
We will denote by $c_{\ell}^{(p)}$ the optimal positive constant such that for every $p\in [1,+\infty)$, $f\in L^p(\O)$ and  $u\in W^{1,p}(\O)$ satisfying
$$
\Delta u=f\quad \text{in }\mathcal{D}'(\O),
$$
there holds
$$
\Vert u\Vert _{W^{2,p}(\O)}\leq c_{\ell}^{(p)}\left(\Vert f\Vert _{L^p(\O)}+\Vert u\Vert _{L^p(\O)}\right).
$$
The optimal constant in the Poincar\'e-Wirtinger inequality will be denoted by $C_{PW}^{(p)}(\O)$. This inequality reads: for every $u\in W^{1,p}(\O)$,
\begin{equation}\label{PW}
\bigg\Vert u -\fint_\O u\bigg\Vert_{L^p(\O)}\leq C_{PW}^{(p)}(\O)\Vert \n u\Vert _{L^p(\O)}.
\end{equation}
We will also use the following regularity results:
\begin{thmnonumbering}(\cite[Theorem 9.1]{Stampacchia})
Let $\O$ be a  $\mathscr C^2$ domain. There exists a constant $C_\O>0$ such that, if $f\in L^\infty(\O)$ and $u\in W^{1,2}(\O)$ solve
\begin{equation}
\left\{
\begin{array}{ll}
-\Delta u=f&\text{ in }\O,
\\\frac{\partial u}{\partial \nu}=0&\text{ on }\partial \O,\end{array}\right.\end{equation}then 
\begin{equation}\label{ControleStampacchia}
\Vert \n u\Vert _{L^1(\O)}\leq C_\O \Vert f\Vert_{L^1(\O)}.
\end{equation}
\end{thmnonumbering}
\begin{thmnonumbering}(\cite[Theorem 1.1]{Mazya})\label{MazyaRegularity} Let $\O$ be a  $\mathscr C^2$ domain. There exists a constant $C_\O>0$ such that, if $f\in L^\infty(\O)$ and $u\in W^{1,2}(\O)$ solve
\begin{equation}
\left\{
\begin{array}{ll}
-\Delta u=f&\text{ in }\O,
\\\frac{\partial u}{\partial \nu}=0&\text{ on }\partial \O,\end{array}\right.\end{equation}then 
\begin{equation}\label{ControleMazya}
\Vert \n u\Vert _{L^\infty(\O)}\leq C_\O \Vert f\Vert _{L^\infty(\O)}.\end{equation}\end{thmnonumbering}
\begin{Remark}
This result is in fact a corollary \cite[Theorem 1.1]{Mazya}. In this article it is proved that the $L^\infty$ norm of the gradient of $f$ is bounded by the Lorentz norm $L^{n,1}(\O)$ of $f$, which is automatically controlled by the $L^\infty(\O)$ norm of $f$.

Note that Stampacchia's orginal result deals with Dirichlet boundary conditions. However, the same duality arguments provide the result for Neumann boundary conditions. 
\end{Remark}
\begin{thmnonumbering}(\cite{LionsMagenes})
Let $r\in (1;+\infty)$. There exists $C_r>0$ such that, if $f\in L^r(\O)$ and if $u\in W^{1,r}(\O)$ be a solution of 
\begin{equation}
\left\{
\begin{array}{ll}
-\Delta u=\operatorname{div}(f)&\text{ in }\O,
\\\frac{\partial u}{\partial \nu}=0&\text{ on }\partial \O,\end{array}\right.\end{equation}
then  there holds
\begin{equation}\label{ControleLionsMagenes}
\Vert \n u\Vert _{L^r(\O)}\leq C_r \Vert f\Vert _{L^r(\O)}.\end{equation}
\end{thmnonumbering}

%%%%%%%%%%%%%%%%%%%%%%%%%%%%%%%%%%%%%%%%%%%

\section{Proofs of the main results}\label{sec:proofMainres}
\subsection{First order optimality conditions for Problem \eqref{TSOM}}
To prove the main results, we first need to state the first order optimality conditions for Problem \eqref{TSOM}.
For that purpose, let us introduce the tangent cone to $\mathcal{M}_{m_0,\kappa}(\O)$ at any point of this set.

\begin{definition} (\cite[chapter 7]{henrot-pierre})
For every $m\in \mathcal M_{m_0,\kappa}(\O)$, the tangent cone to the set $\mathcal{M}_{m_0,\kappa}(\O)$ at $m$, denoted by $\mathcal{T}_{m,\mathcal{M}_{m_0,\kappa}(\O)}$ is the set of functions $h\in \mathrm{L}^\infty(\O)$ such that, for any sequence of positive real numbers $\varepsilon_n$ decreasing to $0$, there exists a sequence of functions $h_n\in \mathrm{L}^\infty(\O)$ converging to $h$ as $n\rightarrow +\infty$, and $m+\varepsilon_nh_n\in\mathcal{M}_{m_0,\kappa}(\O)$ for every $n\in\N$.\label{footnote:cone}
\end{definition}

We will show that, for any $m\in \mathcal{M}_{m_0,\kappa}(\O)$ and any admissible perturbation $h\in \mathcal{T}_{m,\mathcal{M}_{m_0,\kappa}(\O)}$,
the functional ${\mathcal F_\mu}$ is twice G\^ateaux-differentiable at $m$ in direction $h$. To do that, we will show that the solution mapping 
$$
\mathcal S:m\in \mathcal{M}_{m_0,\kappa}(\O)\mapsto \tmm\in L^2(\O),
$$
where $\tmm$ denotes the solution of \eqref{LDE}, is twice G\^ateaux-differentiable. In this view, we provide several $L^2(\O)$ estimates of the solution $\tmm$.
\begin{lemma}\label{Differentiability}(\cite{DING2010688})
The mappping $\mathcal S$ is twice G\^ateaux-differentiable.
\end{lemma}

For the sake of simplicity, we will denote by $\dot{\theta}_{m,\mu}=d\mathcal S(m)[h]$ the G\^ateaux-differential of $\tmm$ at $m$ in direction $h$ and by $\ddot{\theta}_{m,\mu}=d^2\mathcal S(m)[h,h]$ its second order derivative at $m$ in direction $h$.

Elementary computations show that $\dot{\theta}_{m,\mu}$ solves the PDE
\begin{equation}\label{DeriveePrem} \left\{
\begin{array}{ll}
\mu\D \dot{\theta}_{m,\mu} +(m-2\tmm)\dot{\theta}_{m,\mu}=-h\tmm & \text{ in }\O,\\
\frac{\partial \dot{\theta}_{m,\mu}}{\partial \nu}=0 & \text{ on }\partial \O,
\end{array}
\right.\end{equation}
whereas $\ddot{\theta}_{m,\mu}$ solves the PDE
\begin{equation}\label{DeriveeSeconde} \left\{
\begin{array}{ll}
\mu\D \ddot{\theta}_{m,\mu} +\ddot{\theta}_{m,\mu}(m-2\tmm)=-2\left(h\dot{\theta}_{m,\mu}-\dot{\theta}_{m,\mu}^2\right) & \text{ in }\O,\\
\frac{\partial \ddot{\theta}_{m,\mu}}{\partial \nu}=0 & \text{ on }\partial \O,
\end{array}
\right.\end{equation}

It follows that, for all $\mu>0$, the application ${\mathcal F_\mu}$ is G\^ateaux-differentiable with respect to $m$ in direction $h$ and  its G\^ateaux derivative writes 
$$
d{\mathcal F_\mu}(m)[h]=\int_\O \dtm.
$$
Since the expression of $d{\mathcal F_\mu}(m)[h]$ above is not workable, we need to introduce the adjoint state $p_{m,\mu}$ to the equation satisfied by $\dot \theta_{m,\mu}$, i.e \color{black} the solution of the equation
\begin{equation}\label{EquationAdjoint} \left\{
\begin{array}{ll}
\mu \D \pmm+\pmm(m-2\tmm)=1 & \text{ in }\O,\\
\frac{\partial \pmm}{\partial \nu}=0 & \text{ on }\partial \O.
\end{array}
\right.\end{equation}
Note that $p_{m,\mu}$ belongs to $W^{1,2}(\O)$ and is unique, according to the Fredholm alternative. In fact,  we can prove the following regularity results on $p_{m,\mu}$: $p_{m,\mu}\in L^\infty(\O)$, $\Vert p_{m,\mu}\Vert _{L^\infty(\O)}\leq M$, where $M$ is uniform in $m\in \mathcal M_{m_0,\kappa}(\Omega)$ and, for any $p\in [1;+\infty)$, $p_{m,\mu}\in W^{2,p}(\O)$, so that Sobolev embeddings guarantee that $p_{m,\mu}\in \mathscr C^{1,\alpha}(\O)$.

Now, multiplying the main equation of \eqref{EquationAdjoint} by $\dot{\theta}_{m,\mu}$ and integrating two times by parts leads to the expression
$$
d{\mathcal F_\mu}(m)[h]=-\int_\O h\tmm\pmm.
$$ 
Now consider a maximizer $m$. For every perturbation $h$ in the cone $\mathcal{T}_{m,\mathcal{M}_{m_0,\kappa}(\O)}$, there holds $ d{\mathcal F_\mu}(m)[h]\geq 0$.
The analysis of such optimality condition is standard in optimal control theory (see for example \cite{MR1155489}) and leads  to the following result.
\begin{proposition}\label{propoptimCions}
Let us define $\phimm=\tmm\pmm$, where $\tmm$ and $\pmm$ solve respectively equations \eqref{LDE} and \eqref{EquationAdjoint}.
There exists $c\in \R$ such that
$$
\{\varphi_{m,\mu}< c\}=\{m=\kappa\}, \quad \{\varphi_{m,\mu}= c \}=\{0<m<\kappa\},\quad \{\varphi_{m,\mu}> c\}=\{m=0\}.
$$
\end{proposition}

\subsection{Proof of Proposition \ref{SaturationContraintePonctuelle}}
An easy but tedious computation shows that the function $\phimm$ introduced in Proposition \ref{propoptimCions} is $\mathscr C^{1,\alpha}(\O)\cap W^{1,2}(\O)$ function, as a product of two $\mathscr C^{1,\alpha}$ functions and satisfies (in a $W^{1,2}$ weak sense)
\begin{equation}\label{EquationSwitch} \left\{
\begin{array}{ll}
  \mu\D \phimm -2\mu \left\langle \n \phimm,\frac{\n \tmm}{\tmm}\right\rangle+\phimm\left(2\mu\frac{|\n \tmm|^2}{\tmm^2}+2m-3\tmm\right)=\tmm & \text{ in }\O,\\
 \frac{\partial \phimm}{\partial \nu}=0 & \text{ on }\partial \O,
\end{array}
\right.\end{equation}
where $\langle\cdot , \cdot \rangle$ stands for the usual Euclidean inner product.
To prove that $|\{m=0\}|+|\{m=\kappa\}|>0$, we argue by contradiction, by assuming that $|\{m=\kappa\}|=|\{m=0\}|=0$. Therefore, $\phimm= c$ a.e. in $\O$ and, according to \eqref{EquationSwitch}, there holds
$$
c\left(2\mu\frac{|\n \tmm|^2}{\tmm^2}+2m-3\tmm\right)=\tmm
$$
Integrating this identity and using that $\tmm> 0$ in $\O$ and $c\neq 0$, we get
$$2c\left(\mu   \int_\O \frac{|\n \tmm|^2}{\tmm^2}+\int_\O(m-\tmm)\right)=(c+1)\int_\O\tmm.$$
Equation \eqref{metz2135} yields that the left-hand side equals 0, so that one has  $c=-1$. Coming back to the equation satisfied by $\phimm$ leads to
$$m=\tmm-\mu\frac{|\n \tmm|^2}{\tmm^2}.$$
The logistic diffusive equation \eqref{LDE} is then transformed into 
$$\mu \tmm\Delta \tmm-\mu|\n \tmm|^2=0.$$ 
Integrating this equation by parts yields $\int_\O |\n \tmm|^2=0$.
Thus, $\tmm$ is constant, and so is $m$. In other words, $m=m_0$, which, according to \eqref{metz2135} (see Remark \ref{rk:min}) is impossible. The expected result follows.

%%%Preuve du th\'eorÃšme principal
\subsection{Proof of Theorem \ref{TheoremePrincipal}}
The proof of Theorem \ref{TheoremePrincipal} is  based on a careful asymptotic analysis with respect to the diffusivity variable $\mu$. 

Let us first explain the outlines of the proof.

Let us fix $m\in \mathcal{M}_{m_0,\kappa}(\O)$ and $h\in L^\infty(\O)$. In the sequel, the dot or double dot notation $\dot{f}$ or $\ddot{f}$ will respectively denote first and second order G\^ateaux-differential of $f$ at $m$ in direction $h$.  

According to Lemma \ref{Differentiability}, ${\mathcal F_\mu}$ is twice G\^ateaux-differentiable and its second order G\^ateaux-derivative is given by
$$
d^2{\mathcal F_\mu}(m)[h,h]=\int_\O \ddot{\theta}_{m,\mu},
$$
where $\ddot{\theta}_{m,\mu}$ is the second G\^ateaux derivative of $\mathcal S$, defined as the unique solution of \eqref{DeriveeSeconde}. 

Let $m_1$ and $m_2$ be two elements of $\mathcal{M}_{m_0,\kappa}(\O)$ and define
$$
\phi_\mu:[0;1]\ni t\mapsto {\mathcal F_\mu}\Big(tm_2+(1-t)m_1\Big)-t{\mathcal F_\mu}(m_2)-(1-t){\mathcal F_\mu}(m_1).
$$
One has
$$
\frac{d^2\phi_\mu}{dt^2}(t)=\int_\O \ddot\theta_{(1-t)m_1+tm_2,\mu},\quad \text{and}\quad \phi_\mu(0)=\phi_\mu(1)=0,
$$
where $\ddot\theta_{(1-t)m_1+tm_2,\mu}$ must be interpreted as a bilinear form from $L^\infty(\O)$ to $W^{1,2}(\O)$, evaluated two times at the same direction $m_2-m_1$.
Hence, to get the strict convexity of ${\mathcal F_\mu}$, it suffices to show that, whenever $\mu$ is large enough,
$$
\int_\O \ddot{\theta}_{tm_2+(1-t)m_1,\mu}> 0
$$ 
as soon as $m_1\neq m_2$ (in $L^\infty(\O)$) and $t\in (0,1)$, or equivalently that $d^2{\mathcal F_\mu}(m)[h,h]>0$ as soon as $m\in \mathcal{M}_{m_0,\kappa}(\O)$ and $h\in L^\infty(\O)$. Note that since $h=m_2-m_1$, it is possible to assume without loss of generality that $\Vert h\Vert_{L^\infty(\O)}\leq 2\kappa$.
\\
The proof is based on an asymptotic expansion of $\tmm$ into a main term and a reminder one, with respect to the diffusivity $\mu$. It is well-known (see e.g. \cite[Lemma 2.2]{LouMigration}) that one has
\begin{equation}
\tmm\xrightarrow[\mu \to \infty]{W^{1,2}(\O)}m_0.
\end{equation}
%The proof is standard and mainly relies on Sobolev embeddings. 
However, since we are working with resources distributions living in $\mathcal M_{m_0,\kappa}(\O)$, such a convergence property does not allow us to exploit it for deriving optimality properties for Problem \eqref{TSOM}.
%is not enough to carry out an asymptotic analysis. 

For this reason, in what follows, we find a first order term in this asymptotic expansion. To get an insight into the proof's main idea, let us first proceed in a formal way, by looking for a function $\etam$ such that 
$$\tmm\approx m_0+\frac\etam\mu$$ as $\mu \to \infty.$
Plugging this formal expansion in \eqref{LDE} and identifying at order $\frac1\mu$ yields that $\etam$ satisfies
\begin{equation*}
\left\{\begin{array}{ll}
\Delta \etam+m_0(m-m_0)=0&\text{ in }\O\, , 
\\
\frac{\partial \etam}{\partial \nu}=0 &\text{ on }\partial \O.\end{array}\right.
\end{equation*}
To make this equation well-posed, it is convenient to introduce the  function $\hetam$ defined as the unqiue solution to the system
\begin{equation*}
\left\{\begin{array}{ll}
\Delta \hetam+m_0(m-m_0)=0&\text{ in }\O\, , 
\\\frac{\partial \hetam}{\partial \nu}=0 &\text{ on }\partial \O,
\\\fint_\O \hetam=0,&\end{array}\right.\end{equation*}
and to determine a constant $\betam$ such that 
$$\etam=\hetam+\betam.$$
In view of identifying the constant $\betam$, we integrate equation \eqref{LDE} to get 
$$\int_\O \tmm(m-\tmm)=0.$$
which yields, at the order $\frac1\mu$,
$$\betam=\frac1{m_0}\fint_\O \hetam(m-m_0)=\frac1{m_0^2}\fint_\O |\n \hetam|^2.$$
Therefore, one has formally 
$$
\fint_\O \tmm\approx m_0+\frac1\mu\fint_\O |\n \hetam|^2.
$$
As will be proved in Step 1 (paragraph \ref{step1}), the mapping $\mathcal{M}_{m_0,\kappa}(\O)\ni m\mapsto \betam$ is convex so that, at the order $\frac1\mu$, the mapping $\mathcal{M}_{m_0,\kappa}(\O)\ni m\mapsto \fint_\O \tmm$ is convex. We will prove the validity of all the claims above, by taking into account remainder terms in the asymptotic expansion above, to prove that the mapping ${\mathcal F_\mu}:m\mapsto \fint_\O \tmm$ is itself convex whenever $\mu$ is large enough.

\begin{remark}
One could also notice that the quantity $\beta_{1,m}$ arose in the recent paper \cite{HeNi}, where the authors determine the large time behavior of a diffusive Lotka-Volterra competitive system between two populations with growth rates $m_{1}$ and $m_{2}$. If $\beta_{1,m_{1}}>\beta_{1,m_{2}}$, then when $\mu$ is large enough, the solution converges as $t\to +\infty$ to the steady state solution of a scalar equation associated  with the growth rate $m_{1}$. In other words, the species with growth rate $m_{1}$ chases the other one.
In the present article, as a byproduct of our results, we maximize the function $m\mapsto \beta_{1,m}$. This remark implies that this intermediate result might find other applications of its own.
\end{remark}

Let us now formalize rigorously the reasoning outlined above, by considering an expansion of the form
$$\tmm=m_0+\frac\etam \mu+\frac{\mathcal R_{m,\mu}}{\mu^2}.$$
Hence, one has for all $m\in \mathcal{M}_{m_0,\kappa}(\O)$,
$$
d^2{\mathcal F_\mu}(m)[h,h]=\frac{1}{\mu}\int_\O \ddot{\eta}_{1,m}+\frac{1}{\mu^2}\int_\O \ddot{\mathcal{R}}_{m,\mu}
$$
We will show that there holds
\begin{equation}\label{minnea0450}
d^2{\mathcal F_\mu}(m)[h,h]\geq \frac{C(h)}{\mu}\left(1-\frac{\Lambda}{\mu}\right)
\end{equation}
for all $\mu>0$, where $C(h)$ and $\Lambda$ denote some positive constants. 

The strict convexity of ${\mathcal F_\mu}$ will then follow.  Concerning the bang-bang character of maximizers, notice that the admissible set $\mathcal{M}_{m_0,\kappa}(\O)$ is convex, and that its extreme points are exactly the bang-bang functions of $\mathcal{M}_{m_0,\kappa}(\O)$. Once the strict convexity of ${\mathcal F_\mu}$ showed, we then easily infer that ${\mathcal F_\mu}$ reaches its maxima at extreme points, in other words that any maximizer is bang-bang. Indeed, assuming by contradiction the existence of a maximizer writing $tm_1+(1-t)m_2$ with $t\in (0,1)$, $m_1$ and $m_2$, two elements of $\mathcal{M}_{m_0,\kappa}(\O)$ such that $m_1\neq m_2$ on a positive Lebesgue measure set, one has 
$$
{\mathcal F_\mu}(tm_1+(1-t)m_2)< t{\mathcal F_\mu}(m_1)+(1-t){\mathcal F_\mu}(m_2)<\max\{{\mathcal F_\mu}(m_1),{\mathcal F_\mu}(m_2)\},
$$
by convexity of ${\mathcal F_\mu}$, whence the contradiction.

The rest of the proof is devoted to the proof of the inequality \eqref{minnea0450}. It is divided into the following steps:
\begin{itemize}
\item[\bf Step 1.] Uniform estimate of $\int_\O \ddot{\eta}_{1,m}$ with respect to $\mu$.
\item[\bf Step 2.] Definition and expansion of the reminder term $\mathcal{R}_{m,\mu}$.
\item[\bf Step 3.] Uniform estimate of $\mathcal{R}_{m,\mu}$ with respect to $\mu$.
\end{itemize}

\paragraph{Step 1: minoration of $\int_\O \ddot{\eta}_{1,m}$.}\label{step1}
One computes successively
\begin{equation}\label{paris:1827IM}
\dot{\beta}_{1,m}=\frac1{m_0}\fint_\O\left(\dhetam m+\hetam h\right),\quad \ddot{\beta}_{1,m}=\frac1{m_0}\fint_\O \left(2\dhetam h+\ddot{\widehat{\eta}}_{1,m} h\right)
\end{equation}
where $\dhetam$ solves the equation 
\begin{equation}\label{ee}
\left\{
\begin{array}{ll}
  \Delta \dhetam+m_0h=0 & \text{ in }\O\\ 
  \frac{\partial \dhetam}{\partial \nu}=0,& \text{ on }\partial\O
\end{array}
\right. \quad \text{ with }  \int_\O \dhetam=0.
\end{equation}
Notice moreover that $\ddot{\widehat{\eta}}_{1,m}=0$, since $\dhetam$ is linear with respect to $h$. Moreover, multiplying the equation above by $\dhetam$ and integrating  by parts yields
\begin{equation}\label{cvx}
\ddot{\beta}_{1,m}=\frac{2}{m_0^2}\fint_\O |\n \dot{\hat{\eta}}_{1,m}|^2>0
\end{equation}
whenever $h\neq 0$, according to \eqref{paris:1827IM}. Finally, we obtain 
$$
\int_\O \ddot{\eta}_{1,m}=|\O|\ddot{\beta}_{1,m}+\int_\O \ddot{\widehat{\eta}}_{1,m}=|\O|\ddot{\beta}_{1,m}=\frac{2}{m_0^2}\int_\O |\n \dot{\hat{\eta}}_{1,m}|^2.
$$
It is then notable that $\int_\O \ddot{\eta}_{1,m}\geq 0$.

\paragraph{Step 2: expansion of the reminder term $\mathcal{R}_{m,\mu}$.}
Instead of studying directly the equation \eqref{DeriveeSeconde}, our strategy consists in providing a well-chosen expansion of $\ddot{\theta}_{m,\mu}$ of the form 
$$
\ddot{\theta}_{m,\mu}=\sum_{k=0}^{+\infty}\frac{\zeta_k}{\mu^k}, \text{where the $\zeta_k$  are such that }
\sum_{k=2}^{+\infty}\frac{\fint_\O \zeta_k}{\mu^{k-1}}\leq       M\fint_\O \ddot{\eta}_{1,m}.
$$
For that purpose, we will expand formally $\tmm$ as
\begin{equation}\label{DeveloppementSolution}
\tmm=\sum_{k=0}^{+\infty}\frac{\etakm}{\mu^k}.
\end{equation}
Note that, as underlined previously, since $\tmm\underset{\mu\rightarrow +\infty}\longrightarrow  m_0$ in $L^\infty(\O)$, we already know that $\eta_{0,m}=m_0$.

Provided that this expansion makes sense and is (two times) differentiable term by term (what will be checked in the sequel) in the sense of G\^ateaux, we will get the following expansions
$$
\dot{\theta}_{m,\mu}=\sum_{k=0}^{+\infty}\frac{\dot {\eta}_{k,m}}{\mu^k}\quad \text{and}\quad \ddot{\theta}_{m,\mu}=\sum_{k=0}^{+\infty}\frac{\ddot {\eta}_{k,m}}{\mu^k}.
$$
Plugging the expression \eqref{DeveloppementSolution} of $\tmm$ into the logistic diffusive equation  \eqref{LDE}, a formal computation first yields
$$\Delta \etam+m_0(m-m_0)=0, \frac{\partial \eta_{1,m}}{\partial \nu}=0\text{ on }\partial \O$$
$$\Delta \eta_{2,m}+\eta_{1,m}(m-2m_0)=0\text{ in }\O\,, \frac{\partial \eta_{2,m}}{\partial \nu}=0\text{ on }\partial \O$$ and, for any $k\in \N\, , k\geq 2$, $\etakm$ satisfies the induction relation
\begin{equation}\label{eq:etak1init}
\Delta \eta_{k+1,m}+(m-2m_0)\etakm-\sum_{\ell=1}^{k-1}\eta_{\ell,m}\eta_{k-\ell,m}=0\quad \text{ in }\O,
\end{equation}
as well as homogeneous Neumann boundary conditions. These relations do not allow to define $\etakm$ in a unique way (it is determined up to a constant). We introduce the following equations to overcome this difficulty: first, we define $\hat\eta_{1,m}$ and $\hat\eta_{2,m}$ as  the solutions to $$\Delta \hat{\eta}_{1,m}+m_0(m-m_0)=0, \frac{\partial \hat\eta_{1,m}}{\partial \nu}=0\text{ on }\partial \O\, , \fint_\O \hat\eta_{1,m}=0,$$
$$\Delta \hat\eta_{2,m}+\eta_{1,m}(m-2m_0)=0\text{ in }\O\,, \frac{\partial \hat\eta_{2,m}}{\partial \nu}=0\text{ on }\partial \O\, , \fint_\O \hat \eta_{2,m}=0$$ 
and, for any $k\in \N\, , k\geq 2$, we define $\hat{\eta}_{k+1,m}$ as the solution of the PDE
\begin{equation}\label{HierarchySolutionHomogeneous}\left\{
\begin{array}{ll}
 \Delta \hat{\eta}_{k+1,m}+(m-2m_0)\etakm-\sum_{\ell=1}^{k-1}\eta_{\ell,m}\eta_{k-\ell,m}=0 & \text{ in }\O\\ 
 \frac{\partial \hat{\eta}_{k+1,m}}{\partial \nu}=0 & \text{ on }\partial \O
\end{array}
\right.\quad \text{with } \int_\O \hat{\eta}_{k+1,m}=0,
\end{equation}
and to define the real number $\betakm$ in such a way that
\begin{equation}\label{decompetakmhatbeta}
\etakm=\hetakm+\betakm.
\end{equation}
for every $k\in \N^*$. Integrating the main equation of \eqref{LDE} yields 
$$
\int_\O \tmm(m-\tmm)=0.
$$
Plugging the expansion \eqref{DeveloppementSolution} and identifying the terms of order $k$ indicates that we must define $\betakm$ by the induction relation
$$
\begin{array}{ll}
\beta_{1,m}=\frac1{m_0^2}\fint_\O|\n \hat\eta_{1,m}|^2,
\\\beta_{2,m}=\frac1{m_0}\fint_\O m\hat\eta_{2,m}-\frac1{m_0}\fint_\O \eta_{1,m}^2,
\\\beta_{k+1,m}=\frac1{m_0}\fint_\O m\hat\eta_{k+1,m}-\frac1{m_0}\sum_{\ell=1}^{k}\fint_\O \eta_{\ell,m}\eta_{k+1-\ell,m}.\quad (k\geq 2)
\end{array}
$$
This leads to the following cascade system for $\{\hat{\eta}_{k,m},\beta_{k,m}\,,\eta_{k,m}\}_{k\in \N}$:

\begin{equation}\label{HierarchySolution}
{\left\{
\begin{array}{l}
\hat\eta_{0,m}=0,\\ 
\D \hat\eta_{1,m}+m_0(m-m_0)=0\ \text{ in }\O,
\\\D \hat\eta_{2,m}+\eta_{1,m}(m-2m_0)=0\ \text{ in }\O,
\\ \Delta \hat \eta_{k+1,m}+(m-2m_0)\etakm-\sum_{\ell=1}^{k-1}\eta_{\ell,m}\eta_{k-\ell,m}=0\text{ in }\O,\quad (k\geq 2)\\ 
\fint_\O \hat{\eta}_{k,m}=0\,,\quad (k\geq 0)
\\ \frac{\partial \hat{\eta}_{k,m}}{\partial \nu}=0 \text{ over }\partial \O,\quad (k\geq 0)
\\\beta_{0,m}=m_0\, ,
\\\beta_{1,m}=\frac1{m_0^2}\fint_\O |\n \hat \eta_{1,m}|^2\, ,
\\\beta_{2,m}=\frac1{m_0}\fint_\O m\hat\eta_{2,m}-\frac1{m_0}\fint_\O \eta_{1,m}^2,
\\ \beta_{k+1,m}=\frac1{m_0}\fint_\O m\hat\eta_{k+1,m}-\frac1{m_0}\sum_{\ell=1}^{k}\fint_\O \eta_{\ell,m}\eta_{k+1-\ell,m}\, , \quad (k\geq 2)
\\\eta_{k,m}=\hat\eta_{k,m}+\beta_{k,m}.\quad (k\geq 0)           
\end{array}
\right.}
\end{equation}
This implies
$$\fint_\O \eta_{k,m}=\beta_{k,m}\, ,\frac{\partial \eta_{k,m}}{\partial \nu}=0\text{ over }\partial \O.\quad (k\geq 0).$$

Now, the G\^ateaux-differentiability of  both $\hetakm$ and $\betakm$  with respect to $m$ follows from similar arguments as those used to prove Proposition \ref{Differentiability}.
Similarly to System \eqref{HierarchySolution}, the system satisfied by the derivatives needs the introduction of two auxiliary sequences $\{\dot{\hat\eta}_{k,m}\}_{k\in\N}$ and $\{\ddot{\hat\eta}_{k,m}\}_{k\in \N}$. More precisely, we expand $\dot \theta_{m,\mu}$ as 
$$\dot\theta_{m,\mu}=\sum_{k=0}^\infty\frac{\dot \eta_{k,m}}{\mu^k}$$with 
$$\dot\eta_{k,m}=\dot{\hat\eta}_{k,m}+\dot \betakm,$$ 
and  the sequence $\{\dot{\hat\eta}_{k,m}\, , \dot \beta_{k,m}\,, \dot\eta_{k,m}\}_{k\in \N}$  satisfies 
%hi\'erarchie pour la d\'eriv\'ee premiÃšre
\begin{equation}\label{HierarchyDerivative}\left\{
\displaystyle\begin{array}{l}
\dot\eta_{0,m}=0,%\text{ and for all }k\in \N,\\ 
%\Delta \dot\eta_{k+1,m}+(m-2m_0)\detakm-2\sum_{\ell=1}^{k-1}\dot\eta_{\ell,m}\eta_{k-\ell,m}=-h\ekm\text{ in }\O,\\ 
%\\
\\\Delta \dot{\hat \eta}_{1,m}+m_0h=0\text{ in } \O,\\
\Delta\dot{\hat\eta}_{2,m}+\dot\eta_{1,m}(m-2m_0)=-h\dot \eta_{1,m}\text{ in }\O,
\\\Delta \dot{\hat\eta}_{k+1,m}+(m-2m_0)\detakm-2\sum_{\ell=1}^{k-1}\dot\eta_{\ell,m}\eta_{k-\ell,m}=-h\ekm\text{ in }\O,\quad (k\geq 2)\\ 
\fint_\O \dot{\hat\eta}_{k,m}=0, \quad (k\geq 0)
\\\frac{\partial \dot{\hat\eta}_{k,m}}{\partial \nu}=0 \text{ over }\partial \O\,,\quad (k\geq 0)
\\\dot\beta_{0,m}=0\, ,\\
\dot\beta_{1,m}=\frac1{m_0}\fint_\O\left(h\hetam +m\dhetam \right)=\frac2{m_0^2}\fint_\O \langle \n \dot \eta_{1,m},\n \eta_{1,m}\rangle,
\\\dot\beta_{2,m}=\frac1{m_0}\fint_\O (h\hat\eta_{2,m}+m\dot{\hat \eta}_{2,m})-\frac2{m_0}\fint_\O \dot\eta_{1,m}\eta_{1,m},
\\\dot{\beta}_{k+1,m}=\frac1{m_0}\fint_\O (h\hat{\eta}_{k+1,m}+m\dot{\hat{\eta}}_{k+1,m})-\frac2{m_0}\sum_{\ell=1}^{k}\fint_\O \dot\eta_{\ell,m}\eta_{k+1-\ell,m},\quad (k\geq 2)
\\\dot\eta_{k,m}=\dot{\hat \eta}_{k,m}+\dot\beta_{k,m}.\quad (k\geq 0)
\end{array}
\right.
\end{equation}
We note that this implies, for any $k\in \N$,
$$\fint_\O \dot\eta_{k,m}=\dot\beta_{k,m}\, , \frac{\partial \dot \eta_{k,m}}{\partial \nu}=0\text{ on }\partial \O.\quad (k\geq 0)$$
Let us also write the system satisfied by the second order differentials. One gets the following hierarchy for $\{\ddot{\hat\eta}_{k,m}\,,\ddot\beta_{k,m}\,,\ddot\eta_{k,m}\}_{k\in\N}$:
\begin{equation}\label{HierarchySecondDerivative}\left\{
\begin{array}{l}
\ddot{\hat\eta}_{0,m}=0,\\
\Delta \ddot{\hat \eta}_{1,m}=0\text{ in }\O,\\
\Delta \ddot{\hat \eta}_{2,m}+(m-2m_0)\ddot\eta_{1,m}=-2h\dot\eta_{1,m}\text{ in }\O\, ,\\
\Delta \ddot{\hat\eta}_{k+1,m}+(m-2m_0)\ddot\eta_{k,m}-2\sum_{\ell=1}^{k-1}\ddot\eta_{\ell,m}\eta_{k-\ell,m}=2\Big(\sum_{\ell=1}^{k-1}\dot\eta_{\ell,m}\dot\eta_{k-\ell,m}-h\dot{\eta}_{k,m}\Big) \text{ in }\O,\quad (k\geq 2)
\\\fint_\O \ddot{\hat\eta}_{k,m}=0,\quad (k\geq 0)
\\\frac{\partial \ddot{\hat\eta}_{k,m}}{\partial \nu}=0 \text{ on }\partial \O,\quad (k\geq 0)
\\\ddot\beta_{0,m}=0\,,
\\\ddot\beta_{1,m}=\frac{2}{m_0^2}\fint_\O |\n \dot{\hat{\eta}}_{1,m}|^2\, ,
\\\ddot\beta_{2,m}=\frac1{m_0}\fint_\O (2h\dot{\hat\eta}_{2,m}+m\ddot{\hat \eta}_{2,m})-\frac2{m_0}\fint_\O \left((\dot\eta_{1,m})^2+\ddot\eta_{1,m}\eta_{1,m}\right)\,,
\\\ddot{\beta}_{k+1,m}=\frac1{m_0}\fint_\O (2h \dot{\hat{\eta}}_{k+1,m}+m\ddot{\hat{\eta}}_{k+1,m})-\frac2{m_0}\sum_{\ell=1}^{k}\fint_\O (\dot\eta_{\ell,m}\dot\eta_{k+1-\ell,m}+\ddot{\eta}_{\ell,m}\eta_{k+1-\ell,m}),\quad (k\geq 2)
\\\ddot\eta_{k,m}=\ddot{\hat\eta}_{k,m}+\ddot\beta_{k,m}.\quad (k\geq 0)\end{array}
\right.
\end{equation}
This  gives
$$\fint_\O \ddot\eta_{k,m}=\ddot\beta_{k,m}\, ,\frac{\partial \ddot\eta_{k,m}}{\partial \nu}=0\text{ over }\partial \O.\quad (k\geq 0).$$

\paragraph{Step 3: uniform estimates of $\mathcal{R}_{m,\mu}$.}
This section is devoted to proving an estimate on $\ddot{\beta}_{k,m}$, namely
\begin{equation}\label{EstimationCentrale}\left\{
\begin{array}{l}
\forall k\in \N^*\, , \left|\ddot{\beta}_{k,m}\right|\leq \Lambda (k)\ddot{\beta}_{1,m},\\ \text{The power series }  \sum_{k=1}^{+\infty}\Lambda(k)x^k \text{ has a positive convergence radius.}\end{array}
\right.
\end{equation}
This estimate is a key point in our reasoning. Indeed, recall that our goal is to prove that ${\mathcal F_\mu}$ is convex. Assuming that Estimates \eqref{EstimationCentrale} hold true, we expand ${\mathcal F_\mu}$ as follows:
$${\mathcal F_\mu}(m)=\sum_{k=0}^\infty \frac{\fint_\O\etakm}{\mu_k}=\sum_{k=0}^{\infty}\frac{\betakm}{\mu^k}.$$
Differentiating this expression twice with respect to $m$ in direction $h$ yields 
$$\ddot {\mathcal F_\mu}(m)[h,h]=\sum_{k=1}^{\infty}\frac{\ddot\beta_{k,m}}{\mu^k}.$$ Note that the sum starts at $k=1$ since $\beta_{0,m}=m_0$ does not depend on $m$.

We can then write
\begin{align*}
\mu \ddot {\mathcal F_\mu}(m)[h,h]&=\ddot{\beta}_{1,m}+\sum_{k=2}^\infty\frac{\ddot\beta_{k,m}}{\mu^{k-1}}\geq \ddot \beta_{1,m}\left(1-\sum_{k=2}^\infty \frac{\Lambda(k)}{\mu^{k-1}}\right)
%\\&=\ddot F_1(m)[h,h]\left(1-\frac1\mu\sum_{k=2}^\infty \frac{\Lambda(k)}{\mu^{k-2}}\right)
\\&=\ddot{\beta}_{1,m}\left(1-\frac1\mu\sum_{k=0}^\infty \frac{\Lambda(k+2)}{\mu^{k}}\right)
\end{align*}
Recall that $\ddot{\beta}_{1,m}$ is positive as soon as $h$ is not identically equal on 0, according to \eqref{cvx}. The power series associated with $\{\Lambda(k+2)\}_{k\in \N}$ also has a positive convergence radius. Then, the right hand side term is positive provided that $\mu$ be large enough. For the sake of notational clarity, we define $\delta$ as follows: by the Rellich-Kondrachov embedding theorem, see \cite[Theorem 9.16]{Brezis},   there exists $\delta >0$ such that the  continuous embedding $W^{1,2}(\O)\hookrightarrow L^{2+\delta}(\O)$ holds. We fix such a $\delta>0$.
\\

Recall that we know from Equation \eqref{cvx} that $\ddot \beta_{1,m}$ is proportional to $||\n \dot \eta_{1,m}||_{L^2(\O)}^2$. From the explicit expression of $\ddot\beta_{k,m}$ in \eqref{HierarchySecondDerivative}, one claims that \eqref{EstimationCentrale} follows both from the positivity of $\ddot{\beta}_{1,m}$ and from the following estimates: 
\begin{empheq}[left=\empheqlbrace]{align}
\tag{$I_\alpha^k$}\label{EstimationAlpha}  \Vert \ekm\Vert_{L^\infty(\O)}\, , \Vert \n \ekm\Vert _{L^\infty(\O)}\leq \alpha(k),  
\\\tag{$I_\sigma^k$}\label{EstimationSigma}  \Vert \n\dot{\eta}_{k,m}\Vert_{L^2(\O)}\leq \sigma(k)\Vert \n \dot{\eta}_{1,m}\Vert _{L^2(\O)},
\\\tag{$I_\gamma^k$}\label{EstimationGamma} \Vert \dot{\eta}_{k,m}\Vert_{L^2(\O)}\leq \gamma(k) \Vert \n{\dot \eta}_{1,m}\Vert _{L^2(\O)},
\\ \tag{$I_{\tilde\gamma}^k$}\label{EstimationGammaTilde} \Vert \dot{\eta}_{k,m}\Vert_{L^{2+\delta}(\O)}\leq \tilde\gamma(k) \Vert \n{\dot \eta}_{1,m}\Vert _{L^2(\O)},
\\\tag{$I_\delta^k$}\label{EstimationDelta}  \Vert \n\ddot{\eta}_{k,m}\Vert_{L^1(\color{black}\O)}\leq \delta(k)\Vert \n \dot{\eta}_{1,m}\Vert _{L^2(\O)}^2,
\\\tag{$I_\e^k$}\label{EstimationEpsilon}  \Vert \ddot{\eta}_{k,m}\Vert_{L^1(\O)}\leq \e(k)\Vert \n\dot{\eta}_{1,m}\Vert _{L^2(\O)}^2.
\end{empheq}
where for all $k\in \N$, the numbers $\alpha(k)$, $\sigma(k)$, $\gamma(k)$, $\gamma(k)$, $\delta(k)$ and $\varepsilon(k)$ are positive.

In what follows, we will write $f  \lesssim g$ when there exists a constant $C$ (independent of $k$) such that $f\leq Cg$.

The end of the proof is devoted to proving the aforementioned estimates. In what follows, we will mainly deal with the indices $k\geq 3$. Indeed, the case $k=2$ is much simpler since, according to the cascade systems \eqref{eq:etak1init}-\eqref{decompetakmhatbeta}-\eqref{HierarchySolution}-\eqref{HierarchyDerivative}-\eqref{HierarchySecondDerivative}, the equations on $\eta_{k,m}$, $\dot \eta_{k,m}$ and $\ddot \eta_{k,m}$ for $k\geq 3$ involve more terms than the ones on $\eta_{2,m}$, $\dot \eta_{2,m}$ and $\ddot \eta_{2,m}$.

\paragraph{\sc Estimate \eqref{EstimationAlpha}}

This estimate follows from  an iterative procedure. 

Let us fix $\alpha(0)=m_0$ and assume that, for some $k\in \N^*$, the estimate $\eqref{EstimationAlpha}$ holds true.

By $W^{2,p}(\O)$ elliptic regularity theorem, there holds
$$
\Vert \etakkm\Vert _{W^{2,p}(\O)}\ls \Vert \etakkm\Vert _\lp+\bigg\Vert (m-2m_0)\etakm-\sum_{\ell=1}^{k-1}\eta_{\ell,m}\eta_{k-\ell,m}\bigg\Vert_\lp.
$$ 
One  thus gets from the induction hypothesis
$$
\bigg\Vert (m-2m_0)\etakm-\sum_{\ell=1}^{k-1}\eta_{\ell,m}\eta_{k-\ell,m}\bigg\Vert_\lp\ls
 \kappa \alpha(k)+\sum_{\ell=0}^{k-1}\alpha(\ell)\alpha(k-\ell).
$$
Moreover, using that $\Vert \etakkm\Vert _{L^p(\O)}\leq \Vert \hat{\eta}_{k+1,m}\Vert _{L^p(\O)}+|\beta_{k+1,m}|$ and the $L^p$-Poincar\'e-Wirtinger inequality (see Section \ref{sec:toolsNota}), we get
$$\Vert \hat{\eta}_{k+1,m}\Vert _\lp\ls \Vert \n \hat{\eta}_{k+1,m}\Vert _\lp.$$
We now use the  result from \cite[Theorem 1.1]{Mazya} recalled in the introduction: it readily yields
$$
\Vert \n \etakkm\Vert _{L^\infty(\O)} \ls \left|\left| \etakm\right|\right|_{L^\infty(\O)}+\bigg \Vert\sum_{\ell=0}^k \eta_{\ell,m}\eta_{k-\ell,m}\bigg \Vert_{L^\infty(\O)} \ls \sum_{\ell=0}^k\alpha(\ell)\alpha(k-\ell).
$$    
The term $\betakkm$ is controlled similarly, so that
$$
|\betakkm|\ls \sum_{\ell=0}^k \alpha(\ell)\alpha(k-\ell)+\sum_{\ell=1}^k \alpha(\ell)\alpha(k+1-\ell).
$$
Since it is clear that the sequence $\{\alpha(k)\}_{k\in \N}$ can be assumed to be increasing, we write
\begin{align*}
\sum_{\ell=0}^k \alpha(\ell)\alpha(k-\ell)+\sum_{\ell=1}^k \alpha(\ell)\alpha(k+1-\ell)&=\sum_{\ell=0}^{k-1}\alpha(k-\ell)\big(\alpha(\ell)+\alpha(\ell+1)\big)+\alpha(0)\alpha(k)
\\&\ls \sum_{\ell=0}^{k-1}\alpha(\ell+1)\alpha(k-\ell).
\end{align*}
Under this assumption, one has
$$
\Vert \etakkm\Vert _{L^\infty(\O)}\leq |\betakkm|+\Vert \hetakkm\Vert _{L^\infty(\O)}\ls \sum_{\ell=0}^{k}\alpha(\ell+1)\alpha(k-\ell).
$$

This reasoning guarantees the existence of a constant $C_1$, depending only on $\O$, $\kappa$ and $m_0$, such that the sequence defined recursively by $\alpha(0)=m_0$ and
$$
\alpha(k+1)=C_1\sum_{\ell=0}^{k-1}\alpha(\ell+1)\alpha(k-\ell)
$$ 
satisfies the estimate \eqref{EstimationAlpha}.

Setting $a_k=\alpha(k)/C_1^k$ for all $k\in \N$, we know that $\{a_k\}_{k\in \N}$ is a shifted  Catalan  sequence (see   \cite{roman2015introduction}), and therefore, the power series $\sum \alpha(k) x^k$ has a positive convergence radius.

%%%%%Estim\'ee beta et gamma
\paragraph{\sc Estimates \eqref{EstimationSigma} and \eqref{EstimationGamma}.}
Obviously, one can assume that $\sigma(0)=\gamma(0)=0$. One again, we work by induction, by assuming these two estimates known at a given $k\in \N$.  Since \eqref{EstimationSigma} is an estimate on the $\l$-norm of the gradient of $\dot{\eta}_{k+1,m}$, it suffices to deal with $\dhetakkm$. According to the Poincar\'e-Wirtinger inequality, one has $\fint_\O| \dhetakkm|^2\ls \fint_\O |\n \dhetakkm|^2$. Now, using the weak formulation of the equations on $\hetakkm$ and $\dot\eta_{1,m}$, as well as the uniform boundedness of $\Vert h\Vert _\li$, we get
\begin{eqnarray*}
\fint_\O |\n \dot{\hat\eta}_{k+1,m}|^2&=&\fint_\O (m-2m_0)\dot\eta_{k,m}\dot{\hat\eta}_{k+1,m}-2\sum_{\ell=1}^{k-1}\fint_\O\eta_{k-\ell,m}\dot\eta_{\ell,m}\dot{\hat\eta}_{k+1,m}+\fint_\O h \ekm\dot{\hat\eta}_{k+1,m}
\\
&\ls& \Vert \dot\eta_{k,m}\Vert _\l\Vert \dot{\hat\eta}_{k+1,m}\Vert _\l+\sum_{\ell=1}^k \alpha(k-\ell)\Vert \dot{\hat\eta}_{k+1,m}\Vert _\l\Vert \dot\eta_{\ell,m}\Vert _\l+
\\
&& \fint_\O \ekm\langle\n \dot{\eta}_{1,m},\n\dot{\hat\eta}_{k+1,m}\rangle+\fint_\O \dot{\eta}_{k+1,m}\langle \n\dot{\eta}_{1,m},\n  \ekm\rangle
\\
&\ls &\Vert \n \dot{\hat\eta}_{k+1,m}\Vert _\l \Vert \n \dot{\eta}_{1,m}\Vert _\l\Big(\gamma(k)+
\sum_{\ell=1}^k\alpha(k-\ell)\gamma(\ell)+\alpha(k)+\alpha(k)\Big)
\\
&\ls & \Vert \n \dot{\hat\eta}_{k+1,m}\Vert _\l \Vert \n \dot{\eta}_{1,m}\Vert _\l \left(\gamma(k)+\sum_{\ell=0}^k\alpha(k-\ell)\alpha(\ell)\right),
 \end{eqnarray*} 
 where the constants appearing in these inequalities only depend on $\O$, $\kappa$ and $m_0$. It follows that there exists a constant $C_2$ such that, by setting for all $k\in\N$,
 $$
\sigma(k+1)=C_2\left(\gamma(k)+\sum_{\ell=0}^k\alpha(k-\ell)\alpha(\ell)\right),
 $$
the inequality \eqref{EstimationSigma} is satisfied at rank $k+1$.

Let us now state the estimate \eqref{EstimationGamma}. By using the Poincar\'e-Wirtinger inequality, one gets 
\begin{eqnarray*}
\left|\dot{\beta}_{k+1,m}\right|&= & \left|\frac1{m_0}\fint_\O (h\hat{\eta}_{k+1,m}+m\dot{\hat{\eta}}_{k+1,m})-\frac2{m_0}\sum_{\ell=1}^{k}\fint_\O \dot\eta_{\ell,m}\eta_{k+1-\ell,m}\right|
\\
&\ls& \fint_\O \langle \n \dot{\eta}_{1,m},\n\hat{\eta}_{k+1,m}\rangle+\Vert \n \dot{\hat{\eta}}_{k+1,m}\Vert _\l+ \Vert \n \dot{\eta}_{1,m}\Vert _\l\sum_{\ell=1}^k \gamma(\ell)\alpha(k+1-\ell)
\\
&\ls &\alpha(k+1)\Vert \n \etam\Vert _\l+\sigma(k+1)\Vert \n \etam\Vert _\l+\Vert \n \dot{\eta}_{1,m}\Vert _{L^2(\O)}\sum_{\ell=1}^k \gamma(\ell)\alpha(k+1-\ell).
\end{eqnarray*}
Once again, since all the constants appearing in the inequalities depend only on $\O$, $\kappa$ and $m_0$, we infer that one can choose $C_3$ such that, by setting 
$$
\gamma(k+1)=C_3\left(\sigma(k+1)+\alpha(k+1)+\sum_{\ell=1}^k \gamma(\ell)\alpha(k+1-\ell)\right),
$$
the estimate \eqref{EstimationGamma} is satisfied. Notice that, by bounding each term $\alpha(\ell),\ell\leq k$ by $\alpha(k)$ and by using the explicit formula for $\sigma(k+1)$, there exists a constant $C_4$ depending only on $\O$, $\kappa$ and $m_0$ such that 
$$
\displaystyle\gamma(k+1)\leq C_4\sum^k_{\ell=0}\alpha(k+1-\ell)\big( \gamma(\ell)+\alpha(\ell)\big).
$$
Under this form, the same arguments as previously guarantee that the associated power series has a positive convergence radius.

\paragraph{\sc Estimate \eqref{EstimationGammaTilde}}
This is a simple consequence of the Sobolev embedding $W^{1,2}(\O)\hookrightarrow L^{2+\delta}(\O)$. Let $C_\delta>0$ be such that, for any $u\in W^{1,2}(\O)$, 
\begin{equation}\label{eqmetz1140}
\Vert u\Vert _{L^{2+\delta}(\O)}\leq C_\delta \Vert u\Vert _{W^{1,2}(\O)}.
\end{equation}
Then, Estimates \eqref{EstimationSigma} and \eqref{EstimationGamma} rewrite 
$$\Vert\dot\eta_{k,m}\Vert_{W^{1,2}(\O)}\leq\left(\sigma(k)+\gamma(k)\right)\Vert \n\dot\eta_{1,m}\Vert _{L^2(\O)}. $$
and setting
$$
\tilde\gamma(k)= C_\delta \left( \sigma(k)+\gamma(k)\right)
$$ 
concludes the proof of Estimate \eqref{EstimationGammaTilde}.

%%%%Estim\'ee epsilon et delta
\paragraph{\sc Estimates \eqref{EstimationDelta} and \eqref{EstimationEpsilon}.}
For the sake of clarity, let us recall that $k\in \N$ being fixed, according to Systems \eqref{HierarchyDerivative} and \eqref{HierarchySecondDerivative}, the functions $\dot\eta_{kw,m}$ and $\ddot\eta_{k,m}$ satisfy respectively 
$$
\Delta \dot\eta_{k+1,m}+(m-2m_0)\detakm-2\sum_{\ell=1}^{k-1}\dot\eta_{\ell,m}\eta_{k-\ell,m}=-h\ekm\text{ in }\O
$$
and
$$
\Delta \ddot\eta_{k+1,m}+(m-2m_0)\ddot\eta_{k,m}-2\sum_{\ell=1}^{k-1}\ddot\eta_{\ell,m}\eta_{k-\ell,m}=2\Big(\sum_{\ell=1}^{k-1}\dot\eta_{\ell,m}\dot\eta_{k-\ell,m}-h\dot{\eta}_{k,m}\Big) \text{ in }\O
$$

 As previously, we first set $\delta(0)=\varepsilon(0)=0$ and argue by induction. 

To prove these estimates, let us first control $\Vert \n \ddot {\hat \eta}_{k+1,m}\Vert _{L^1(\O)}$. To this aim, let us use Estimates \eqref{EstimationGammaTilde}, the Stampacchia regularity Estimate \eqref{ControleStampacchia} and the Lions-Magenes regularity Estimate \eqref{ControleLionsMagenes}. We first use the equation 
$$\Delta \dot \eta_{1,m}+m_0h=0$$ to split the equation on $\ddot{\hat\eta}_{k+1,m}$ in System \eqref{HierarchySecondDerivative} as follows:
$$h\dot \eta_{k,m}=-\frac1{m_0}\dot\eta_{k,m}\Delta \dot \eta_{1,m}=-\frac1{m_0}\left( \text{ div}\left(\dot \eta_{k,m}\n \dot\eta_{1,m}\right)-\langle \dot \n \eta_{k,m},\dot\n \eta_{1,m}\rangle\right).$$
Introduce the function
$$
H_k=(m-2m_0)\ddot\eta_{k,m}-2\sum_{\ell=1}^{k-1}\ddot\eta_{\ell,m}\eta_{k-\ell,m}-2\sum_{\ell=1}^{k-1}\dot\eta_{\ell,m}\dot\eta_{k-\ell,m}+\frac2{m_0}\langle\dot\n \eta_{k,m},\dot\n \eta_{1,m}\rangle
$$
then $\ddot{ \hat\eta}_{k+1,m}$ solves
$$\Delta \ddot{ \hat \eta}_{k+1,m}+H_k=\frac{2}{m_0}\text{ div}\left(\dot\eta_{k,m}\n \dot \eta_{1,m}\right),$$
along with Neumann boundary conditions, according to \eqref{EstimationCentrale}.

By using the induction assumption and the Cauchy-Schwarz inequality, one gets 
$$
\Vert H_k\Vert _{L^1(\O)}\ls \left( \e(k)+\sum_{\ell=1}^{k-1}\e(\ell)\alpha(k-\ell)+\sum_{\ell=1}^{k-1}\gamma(\ell)\gamma(k-\ell)+\gamma(1)\gamma(k)\right)\Vert \n \dot \eta_{1,m}\Vert _{L^2(\O)}^2.
$$
Furthermore, let us consider the same number $\delta>0$ as the one introduced and used in Estimate \eqref{EstimationGammaTilde}, and define $r>1$ such that $\frac1r=\frac12+\frac1{2+\delta}$, where $\delta>0$ is fixed so that \eqref{eqmetz1140} holds true. 
 By combining Estimate \eqref{EstimationGammaTilde} with the H\"{o}lder's inequality, we have
\begin{equation}\label{EstimationLr}\Vert \dot \eta_{k,m}\n \dot \eta_{1,m}\Vert _{L^{r}}\leq \Vert \dot\eta_{k,m}\Vert _{L^{2+\delta}(\O)}\Vert \n \dot \eta_{1,m}\Vert _{L^2(\O)}\leq \tilde\gamma(k)\Vert \n \dot \eta_{1,m}\Vert _{L^2(\O)}^2. \end{equation}

Let us introduce $(\psi_{k+1},\xi_{k+1})$ as the respective solutions of
\begin{equation}\label{EquationPsi} \left\{
\begin{array}{ll}
\Delta \psi_{k+1}+H_k=0 & \text{ in }\O,\\
\frac{\partial \psi_{k+1}}{\partial \nu}=0 & \text{ on }\partial \O,
\\\fint_\O \psi_{k+1}=0,
\end{array}
\right.\end{equation}
and
\begin{equation}\label{EquationXi} \left\{
\begin{array}{ll}
\Delta \xi_{k+1}=-2\text{div}(\dot\eta_{k,m}\n \dot\eta_{1,m}) & \text{ in }\O,\\
\frac{\partial \xi_{k+1}}{\partial \nu}=0 & \text{ on }\partial \O,
\\\fint_\O \xi_{k+1}=0,
\end{array}
\right.\end{equation}
so that $\ddot{\hat\eta}_{k+1,m}=\psi_{k+1}+\xi_{k+1}$. 
Stampacchia's Estimate \eqref{ControleStampacchia} leads to
$$
\Vert \n \psi_{k+1}\Vert _{L^1(\O)}\ls \Vert H_k\Vert _{L^1(\O)}\ls \left( \e(k)+\sum_{\ell=1}^{k-1}\e(\ell)\alpha(k-\ell)+\sum_{\ell=1}^{k-1}\gamma(\ell)\gamma(k-\ell)+\gamma(1)\gamma(k)\right)\Vert \n \dot \eta_{1,m}\Vert _{L^2(\O)}^2 ,
$$
and moreover,
\begin{align*}
\Vert \n \xi_{k+1}\Vert _{L^1(\O)}&{\ls} \Vert \n \xi_{k+1}\Vert _{L^r(\O)}\text{ by H\"{o}lder's inequality}
\\&\ls \Vert \dot \eta_{k,m}\n \dot \eta_{1,m}\Vert _{L^r(\O)}\text{ by Lions and Magenes Estimate }\eqref{ControleLionsMagenes}
\\&\ls \tilde \gamma(k)\Vert \n \dot \eta_{1,m}\Vert _{L^2(\O)}^2 \text{ by Estimate \eqref{EstimationLr}.} \end{align*}
We then have
\begin{multline*}\Vert \n \ddot{\hat\eta}_{k+1,m}\Vert _{L^1(\O)}=\Vert \n \psi_{k+1}+\n \xi_{k+1}\Vert _{L^1(\O)}\\
\ls \left( \e(k)+\sum_{\ell=1}^{k-1}\e(\ell)\alpha(k-\ell)+\sum_{\ell=1}^{k-1}\gamma(\ell)\gamma(k-\ell)+\gamma(1)\gamma(k)+\tilde \gamma(k)\right)\Vert \n \dot \eta_{1,m}\Vert _{L^2(\O)}^2.\end{multline*}
and we conclude by setting $\delta(k+1)=\e(k)+\sum_{\ell=1}^{k-1}\e(\ell)\alpha(k-\ell)+\sum_{\ell=1}^{k-1}\gamma(\ell)\gamma(k-\ell)+\gamma(1)\gamma(k)+\tilde \gamma(k)$.

\medskip

Let us now derive $\e(k+1)$. We proceed similarly to the proof of Estimate \eqref{EstimationGamma}: from the Poincar\'e-Wirtinger Inequality, there holds
$$\left|\left| \ddot \eta_{k+1,m}-\fint_\O \ddot \eta_{k+1,m}\right|\right|_{L^2(\O)}\ls ||\n \ddot \eta_{k+1,m}||_{L^2(\O)}$$ so that, from Estimate \eqref{EstimationDelta} it suffices to control $\fint_\O \ddot \eta_{k+1,m}$.
\\Starting from the expression
$$
\fint_\O \ddot{\eta}_{k+1,m}=\frac1{m_0}\fint_\O (2h \dot{\hat{\eta}}_{k+1,m}+m\ddot{\hat{\eta}}_{k+1,m})-\frac2{m_0}\sum_{\ell=1}^{k}\fint_\O (\dot\eta_{\ell,m}\dot\eta_{k+1-\ell,m}+\ddot{\eta}_{\ell,m}\eta_{k+1-\ell,m})
$$
stated in \eqref{HierarchySecondDerivative} and using the Cauchy-Schwarz inequality, one gets 
\begin{eqnarray*}
\left|\fint_\O \ddot{\eta}_{k+1,m}\right|&=&\left|\frac1{m_0}\fint_\O (2h \dot{\hat{\eta}}_{k+1,m}+m\ddot{\hat{\eta}}_{k+1,m})-\frac2{m_0}\sum_{\ell=1}^{k}\fint_\O (\dot\eta_{\ell,m}\dot\eta_{k+1-\ell,m}+\ddot{\eta}_{\ell,m}\eta_{k+1-\ell,m})\right|\\
&\ls &\left| \fint_\O h \dot{\hat{\eta}}_{k+1,m}\right|+\Vert \ddot{\hat{\eta}}_{k+1,m}\Vert_{L^2(\O)}+\sum_{\ell=1}^{k} \Vert\dot\eta_{\ell,m}\Vert_{l^2(\O)}\Vert\dot\eta_{k+1-\ell,m}\Vert_{L^2(\O)}\\
&& +\sum_{\ell=1}^{k} \Vert \ddot{\eta}_{\ell,m}\Vert_{L^2(\O)}\Vert \eta_{k+1-\ell,m}\Vert_{L^2(\O)}
\end{eqnarray*}
We then use Equation \eqref{ee} to get
$$
\fint_\O h \dot{\hat \eta}_{k+1,m}=\frac1{m_0}\fint_\O\langle \n \dot\eta_{1,m},\n \dot{\hat \eta}_{k+1,m}\rangle\leq \frac1{m_0}\sigma(1)\sigma(k+1)\Vert \n \dot\eta_{1,m}\Vert _\l^2.
$$
Since $\ddot \beta_{1,m}$ is proportional to $||\n \dot \eta_{1,m}||_{L^2(\O)}^2$, this gives
$$\left|\fint_\O \ddot \eta_{k+1,m}\right|\ls\left(\sigma(1)\sigma(k+1)+\delta(k+1)+\sum_{\ell=1}^k \left(\gamma(\ell)\gamma(k+1-\ell)+\alpha(k+1-\ell)\e(\ell)\right)\right)\ddot\beta_{1,m}.
$$
Setting $\e(k+1)=\sigma(1)\sigma(k+1)+\delta(k+1)+\sum_{\ell=1}^k \left(\gamma(\ell)\gamma(k+1-\ell)+\alpha(k+1-\ell)\e(\ell)\right)$ concludes the proof.

\paragraph{Summary.} We have proved here that the functional ${\mathcal F_\mu}$ has an asymptotic expansion of the form
$${\mathcal F_\mu}(m)=m_0+\frac{\beta_{1,m}}\mu+R_\mu(m),$$
where $m\mapsto \beta_{1,m}=\frac1{m_0^2}\int_\O |\n \eta_{1,m}|^2$ is a strictly convex functional, and where $R_\mu$ satisfies the two following conditions:
\begin{enumerate}
\item $R_\mu=\underset{\mu \to\infty}O\left(\frac1{\mu^2}\right)$ uniformly in $\mathcal M_{m_0,\kappa}(\O)$,
\item $R_\mu$ can be expanded in a power series of $\frac1\mu$ as follows:
$$R_\mu(m)=\frac1\mu\sum_{k=2}^\infty\frac{\beta_{k,m}}{\mu^{k-1}}$$,
\item $R_\mu$ is twice G\^ateaux-differentiable, and, for any $m\in \mathcal M_{m_0,\kappa}(\O)$, for any admissible variation $h\in \mathcal T_{m,\mathcal M_{m_0,\kappa}(\O)}$,
$$\mu\left|\ddot R_\mu [h,h]\right|\ls\,   \ddot \beta_{1,m}.$$\end{enumerate}
It immediately follows that the functional ${\mathcal F_\mu}$ satisfies the following lower bound on its second derivative: for any $m\in \mathcal M_{m_0,\kappa}(\O)$, for any admissible variation $h\in \mathcal T_{m,\mathcal M_{m_0,\kappa}(\O)}$,
\begin{equation}\left(1-\frac1\mu\right) \ddot \beta_{1,m}\ls \text{ }  \mu\ddot {\mathcal F_\mu}(m)[h,h],\end{equation}
so that it has a positive second derivative, according to \eqref{cvx}. Hence, ${\mathcal F_\mu}$ is strictly convex for $\mu$ large enough.
\\Since the maximizers of a strictly convex functional  defined on a convex set are extreme points, and that the extreme points of $\mathcal M_{m_0,\kappa}(\O)$ are bang-bang functions, this ensures that all maximizers of ${\mathcal F_\mu}$ are bang-bang functions.

%%%ANALYSE DU PRMEIER PROBLEME LIMITE	
\subsection{Proof of Theorem \ref{Concentration}}%%%les cr\'eneaux sont des maximiseurs locaux pour des pertubations L^1
In what follows, it will be convenient to introduce the functional 
$$
F_1:m\mapsto \beta_{1,m}=\frac1{m_0^2}\fint_\O|\n  \eta_{1,m}|^2=\fint_\O \eta_{1,m}
$$
where $\eta_{1,m}$ is defined as a solution to System \eqref{HierarchySolution}. The index in the notation $F_1$ underlines the fact that $F_1$ involves the solution $\eta_{1,m}$.

According to the proof of Theorem \ref{TheoremePrincipal} (Step 1), we already know that $F_1$ is a convex functional on $\mathcal M_{m_0,\kappa}(\O)$.
\subsubsection{Proof of $\Gamma$-convergence property for general domains}
To prove this theorem, we proceed into three steps: we first prove weak convergence, then show that maximizers of the functional $F_1$ are necessarily extreme points of $\mathcal M_{m_0,\kappa}(\O)$ and finally recast $F_1$ using the energy functional $\mathcal E_m$.  Since weak convergence to an extreme point entails strong convergence, this will conclude the proof of the $\Gamma$-convergence property.
\paragraph{Convergence of maximizers.}
For $\mu>0$, let $m_\mu$ be a solution to \eqref{TSOM}. According to Theorem \ref{TheoremePrincipal}, there exists $\mu^*>0$ such that $m_\mu=\kappa \chi_{E_\mu}$ for all $\mu\geq \mu^*$, where $E_\mu\subset \O$ is such that $|E_\mu|=m_0\frac{|\O|}\mu$.
\\Since the family $\{m_\mu\}_{\mu>0}$ is uniformly bounded in $L^\infty(\O)$, it converges up to a subsequence to some element $m_\infty\in \mathcal{M}_{m_0,\kappa}(\O)$,  weakly star in $L^\infty(\O)$.
Observe that the maximizers of ${\mathcal F_\mu}$ over $\mathcal{M}_{m_0,\kappa}(\O)$ are the same as the maximizers of $\mu({\mathcal F_\mu}-m_0)$. Recall that, given $m$ in $\mathcal{M}_{m_0,\kappa}(\O)$, there holds $\mu({\mathcal F_\mu}(m)-m_0)=\fint_\O \eta_{1,m}+ \underset{\mu\to \infty}{\operatorname{O}}(\frac1\mu)$ according to the proof of Theorem \ref{TheoremePrincipal}, where the notation $\operatorname{O}\left(\frac1\mu\right)$ stands for a function uniformly bounded in $L^\infty(\O)$. In other words, we have 
$$\mu\left(\mathcal F_\mu-m_0\right)=F_1+\underset{\mu\to \infty}{\operatorname{O}} \left(\frac1\mu\right)$$ with the same notation for  $\operatorname{O}\left(\frac1\mu\right)$.

For an arbitrary $m\in \mathcal{M}_{m_0,\kappa}(\O)$, by passing to the limit in the inequality  
$$
\mu({\mathcal F_\mu}(m_\mu)-m_0)\geq \mu({\mathcal F_\mu}(m)-m_0)
$$ 
one gets that $m_\infty$ is necessarily a maximizer of the functional
$F_1$ over $ \mathcal{M}_{m_0,\kappa}(\O)$. 

%\paragraph{Step 2: strong convergence of maximizers.}
Wa have shown in the proof of Theorem \ref{TheoremePrincipal} that $F_1$ is convex on $\mathcal M_{m_0,\kappa}(\O)$ (Step 1). Its maximizers are thus extreme points. It follows that any weak limit of $\{m_\mu\}_{\mu>0}$ is an extreme point to this set. Thus, the convergence is in fact strong in $L^1$ (\cite[Proposition 2.2.1]{henrot-pierre}). 

\paragraph{``Energetic'' expression of $F_1(m)$.} 
%We prove here that we can rewrite $F_1$ as an energy functional.
Recall that $F_1$ is given by 
$$F_1(m)=\frac1{m_0^2}\fint_\O |\n \etam|^2,$$
where $\etam$ solves
$$
\left\{\begin{array}{ll}
\Delta \etam+m_0(m-m_0)=0&\text{ in }\O,
\\\frac{\partial \etam}{\partial \nu}=0&\text{ on } \partial \O ,
\\ \fint_\O \etam=\frac1{m_0^2}\fint_\O |\n \etam|^2.&\end{array}\right.$$
The last constraint, which is derived from the integration of Equation \eqref{LDE}, by passing to the limit as $\mu \to +\infty$, is not so easy to handle. This is why we prefer to deal with $\hat\eta_{1,m}$, solving the same equation as $\etam$ completed with the integral condition 
$$
\fint_\O \hat \eta_{1,m}=0.
$$
Since $\etam$ and $\hat \eta_{1,m}$ only differ up to an additive constant, we have $\n \etam=\n \hat\eta_{1,m}$, so that
\begin{equation}
F_1(m)=\frac1{m_0^2}\fint_\O |\n \hat \eta_{1,m}|^2\quad \text{and}\quad \hat \eta_{1,m}\in X.
\end{equation}
Regarding then the variational problem
\begin{equation}
\label{ProblemeLimite}\tag{$PV_1$}\sup_{m\in \mathcal M(\O)}F_1(m),
\end{equation}
and standard reasoning
%Recall that we have introduced the function $$\mathcal E_m:X\ni u\mapsto \frac12\fint_\O |\n u|^2-m_0\fint_\O mu.$$  
on the PDE solved by $\hat \eta_{1,m}$ yields that 
$$F_1(m)=-2\min_{u\in X}\mathcal E_m(u),$$ 
%We can then write $$\max_{m\in \mathcal M_{m_0,\kappa}(\O)} F_1(m)=-2\min_{m\in \mathcal M_{m_0,\kappa}(\O)}\min_{u\in X}\mathcal E_m(u).$$
leading to the desired result.

\subsubsection{Properties of maximizers of $F_1$ in a two-dimensional orthotope}
We investigate here the case of the two-dimensional orthotope $\O=(0;a_1)\times (0;a_2)$.
In the last section, we proved that every maximizer $m$ of $F_1$ over $\mathcal M_{m_0,\kappa}(\O)$ is of the form $m=\kappa \chi_E$ where $E$ is a measurable subset of $\O$ such that $\kappa |E|=m_0|\O|$.

Let $E^*$ be such a set. We will prove that $E^*$ is, up to a rotation of $\O$, decreasing in every direction. It relies on the combination of symmetric decreasing rearrangements properties and optimality conditions for Problem \eqref{ProblemeLimite}. 

Introduce the notation $\hat \eta_{1,E^*}:=\hat \eta_{1,\kappa \chi_{E^*}}$. A similar reasoning to the one used in Proposition \ref{propoptimCions} (see e.g. \cite{MR1155489}) yields the existence of a Lagrange multiplier $c$ such that 
\begin{equation}\label{cionsOptPV1}
\{\hat\eta_{1,E^*}> c\}=E^*,\quad \{\hat\eta_{1,E^*}<c\}=(E^*)^c,\quad \{\hat\eta_{1,E^*}= c\}=\partial E^*.
\end{equation}
We already know, thanks to the equality case in the decreasing rearrangement inequality, that any maximizer $E^*$ is decreasing or increasing in every direction.

To conclude, it remains to prove that $E^*$ is connected. Let us argue by contradiction, by assuming that $E^*$ has at least two connected components. 
%We use two different results to derive a contradiction: first of all, we use the case of equality in the Polya-Szego inequality (see Section \ref{Ptes}) and the optimality conditions on \eqref{ProblemeLimite} to prove that $E^*$ has at most two components. We then use the optimality conditions once again on the variational problem \eqref{ProblemeLimite} to prove that $E^*$ only has one component.

In what follows, if $E$ denotes a measurable subset of $\Omega$, we will use the notation $\hat \eta_{1,E}:=\hat \eta_{1,\kappa \chi_{E}}$. The steps of the proof are illustrated on Figure \ref{figProofPS} below.

\paragraph{Step 1: $E^*$ has at most two components.}
It is clear from the equality case in the P\`olya-Szeg\"o inequality that $\hat \eta_{1,E^*}$ is decreasing in every direction (i.e, it is either nondecreasing or  nonincreasing on every horizontal or vertical line).
%The equality case in the Polya-Szego inequality \\

Let $e_1=(0,0)\, , e_2=(a_1,0)\, , e_3=(a_1,a_2)\, , e_4=(0,a_2)$ be the four vertices of the orthotope $\O=(0;a_1)\times (0;a_2)$.
Let $E_1$ be a connected component of $E^*$. Since $E_1$ is monotonic in both directions $x$ and $y$, thus it necessarily contains at least one vertex. Up to a rotation, one can assume that $e_1\in E_1$.
Since $E_1$ is decreasing in the direction $y$, there exists $\underline x  \in [0;a_1]$ and a non-increasing function $f:[0;\underline x]\to[0;a_1]$ such that
$$
E_1=\left\{(x,t)\, , x\in (0;x_1)\, , t \in [0;f(x)]\right\}
$$
Since $f$ is decreasing, one has $E_1\subseteq [0;\underline x]\times [0;f(0)]$. 

Let $E_2$ be another connected component of $E^*$. Since $E^*$ is monotonic in every direction, the only possibility is that $E_2$ meet the upper corner $[\underline x;a_1]\times [f(0);a_2]$, meaning that $e_3\in E_2$ and therefore, there exist $\overline x\in [\underline x;a_1]$ and a non-decreasing function $g:[\overline x;a_2]\to[0;a_2]$ such that
$$E_2=\left\{(x,t)\, , x\in [\overline x;a_1]\, , t\in [a_2-g(x);a_2]\right\}$$
 %Furthermore, we also deduce that 
%$$\underset{[0;\underline x]}\sup f\geq  \underset{[\overline x;a_1]}\sup g.$$

\paragraph{Step 2: geometrical properties of $E_1$ and $E_2$.} 
We are  going to prove that   $g$ or $f$ is constant  and that $\underline x=\overline x$.
Let $b_2$ be the decreasing rearrangement in the direction $y$. Let $E_*:=b_2(E^*)$. 

We claim that, by optimality of $E^*$, we have
 \begin{equation}\label{SolutionRearrangementt}
 F_1(b_2(E^*))=F_1(E^*)\text{ and } b_2(\hat \eta_{1,E^*})=\hat\eta_{1,b_2(E^*)}.
 \end{equation}
 For the sake of clarity, the proof of \eqref{SolutionRearrangementt} is postponed to the end of this step.
 
Since $b_2(E^*)$ is necessarily a solution of Problem  \eqref{ProblemeLimite}, it follows, by monotonicity of maximizers, that the mapping  $\tilde f:x\in [0;a_1]\mapsto \mathcal H^1\Big(\left(\{x\}\times [0;a_2]\right)\cap b_2(E^*)\Big)$ is also monotonic.
However, it is straightforward that $\tilde f=f\chi_{[0;\underline x]}+g\chi_{[\overline x;a_1]}$.
% and we already know that $f$ is non-increasing, and $g$ is non-decreasing.  
If $f$ is nonconstant, it follows that $\tilde f$ is non-increasing. Since $g$ is non-decreasing and has the same monotonicity as $\tilde f$, it follows that $g$ is necessarily constant. Hence, we get that $\underline x=\overline x$ and that $\underset{[0;\underline x]}\inf f$ is positive. Else, $\tilde f$ would be non-increasing and vanish in $(\underline x;\overline x)$. Finally, we also conclude that  $\underset{[0;\underline x]}\inf f \geq g.$ Thus, we can consider the following situation: $\underline x=\overline x$, $f\geq \alpha$ and $f$ is non-increasing and $g$ is constant, i.e $g=\alpha$.

 \begin{proof}[Proof of \eqref{SolutionRearrangementt}]
 Recall that for every $m\in \mathcal M(\O)$, $\hetam$ is the unique minimizer of the energy functional $\mathcal E_m$ over $X$ where $\mathcal E_m$ and $X$ are defined by \eqref{EspaceNul}-\eqref{def:Em}. For a measurable subset $E$ of $\O$, introduce the notations $F_1(E):=F_1(\kappa \chi_E)$ and $\mathcal E_E:=\mathcal E_{\kappa \chi_E}$. Since $F_1(m)=-\frac2{m_0^2}\mathcal E_m(\hetam)$ and since $E^*$ is a maximizer of $F_1$, we have
 $$F_1(E^*)\geq F_1(b_2(E^*)).$$
Furthermore, one has
 \begin{align*}
 F_1(E^*)&=-\frac2{m_0^2}\mathcal E_{E^*}(\hat\eta_{1,E^*}) \leq -\frac2{m_0^2}\mathcal E_{b_2(E^*)}(b_2(\hat\eta_{1,E^*}))\\
&  \leq -\frac2{m_0^2}\mathcal E_{b_2(E^*)}(\hat \eta_{1,b_2(E^*)})= F_1({b_2(E^*)}).
 \end{align*}
 by using successively the Hardy-Littlewood and P\`olya-Szeg\"o inequalities.
 
 Thus, all these inequalities are in fact equality, which implies that $b_2(E^*)$ is also a maximizer of $F_1$ over $ \mathcal M(\O)$. Furthermore, by the equimeasurability property, one has $b_2(\hat\eta_{1,E^*})\in X$, so that $b_2(\hat\eta_{1,E^*})$ is a minimizer of $\mathcal E_{b_2(E^*)}$ over $X$. The conclusion follows.
 \end{proof}

\paragraph{Step 3: $E^*$ has at most one component.}
To get a contradiction, let us use the optimality conditions \eqref{cionsOptPV1}. This step is illustrated on the bottom of Figure \ref{figProofPS}. 

 By using the aforementioned properties of maximizers, we get that $\eta_{1,b_2(E^*)}$ is constant and equal to $c$ on $\{\underline x\}\times [0;\alpha ]\subset b_2(E^*)$:
 \begin{equation}\label{Minimum}
 \hat \eta_{1,b_2(E^*)}=c \text{ on } \{\underline x\}\times [0;\alpha].
 \end{equation} Furthermore, since $b_2(E^*)$ is a maximizer of $F_1$, it follows that $\hat\eta_{1,b_2(E^*)}$ is constant on $\partial b_2(E^*)$. But one has $b_2(\hat \eta_{1,E^*})=\hat\eta_{1,E^*}$ on $[0;\underline x]\times [0;a_2]$ since $\hat \eta_{1,E^*}$ is decreasing in the vertical direction on this subset.  We get that $\hat\eta_{1,b_2(E^*)}$ is equal to $c$ on $\partial b_2(E^*)$
 
 However, by the strict maximum principle, $\hat \eta_{1,b_2(E^*)}$ cannot reach its minimum in $b_2(E^*)$, which is a contradiction with \eqref{Minimum}. This concludes the proof. 
 %Thus, $E^*$ is connex, and is as a consequence decreasing in every direction (up to a rotation).

\begin{figure}[h]
\begin{center}
\includegraphics[width=5.2cm]{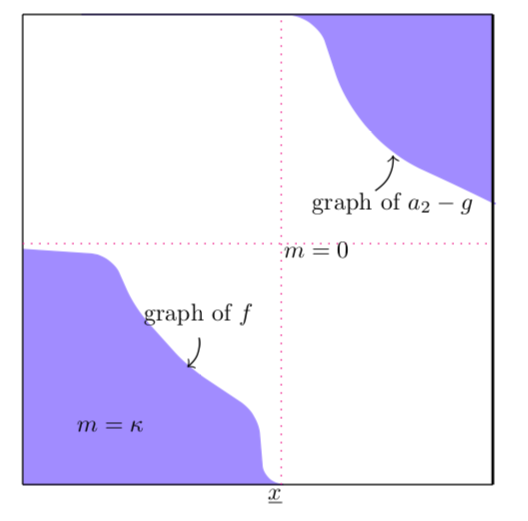}\hspace{1cm}
\includegraphics[width=5.2cm]{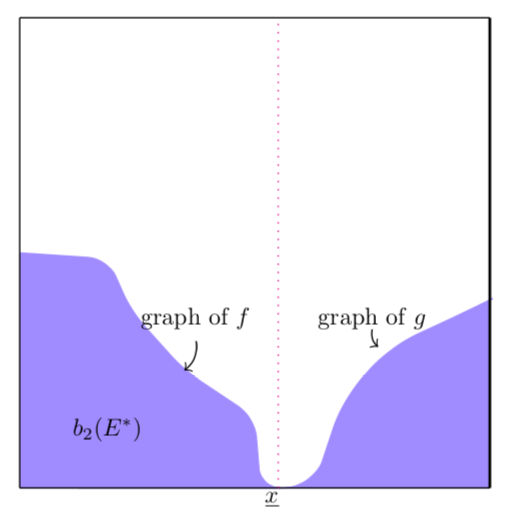}\\
\includegraphics[width=5.2cm]{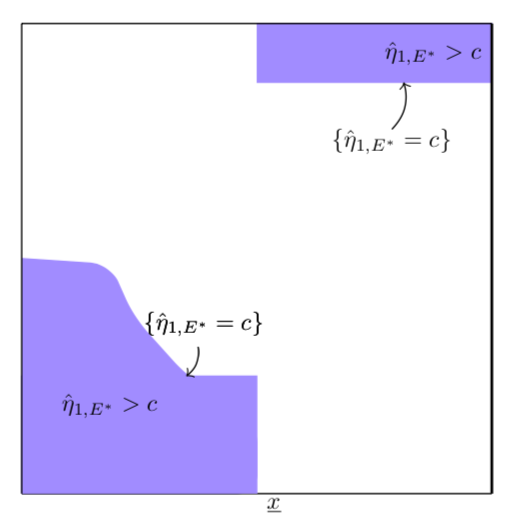}\hspace{1cm}
\includegraphics[width=5.2cm]{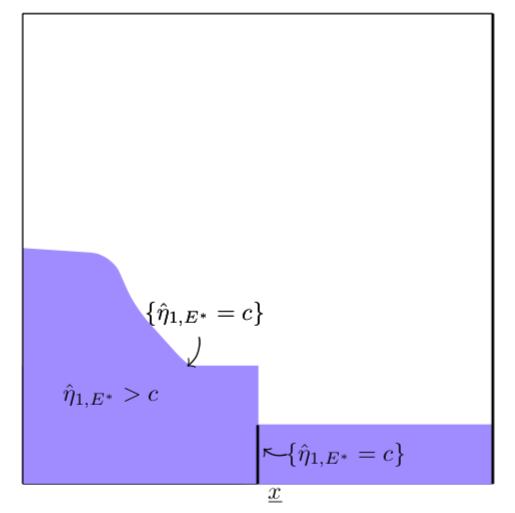}
\caption{Illustration of the proof and the notations used. \label{figProofPS}}
\end{center}
\end{figure}

\subsection{Proof of Theorem \ref{Creneau}}\label{Derniere}
%We will prove that, if $\mu$ large enough, any maximizer of ${\mathcal F_\mu}$ is equal to either $\tilde m=\kappa\chi_{(1-\ell,1)]}$ or $\tilde m(1-\cdot)=\kappa \chi_{(0,\ell)}$ with $\ell=\frac{m_0}\kappa$. 
As a preliminary remark, we claim that the function $\theta_{\tilde m,\mu}$ solving \eqref{LDE} with $m=\tilde m$ is positive increasing. Indeed, recall that $\theta_{\tilde m,\mu}$ is the unique minimizer of the energy functional
\begin{equation}\label{Energie}
\mathcal{E}:W^{1,2}(\O,\R_+)\ni u\mapsto \frac{\mu}{2}\int_0^1u'^2-\frac{1}{2}\int_0^1m^*u^2+\frac{1}{3}\int_0^1u^3.
\end{equation}
By using the rearrangement inequalities recalled in Section \ref{sec:toolsNota} and the relation $(\tilde m)_{br}=\tilde m$\color{black}, one easily shows   that 
$$
\mathcal{E}(\theta_{\tilde m,\mu})\geq \mathcal{E}((\theta_{\tilde m,\mu})_{br}),
$$
and therefore, one has necessarily $\theta_{\tilde m,\mu}=(\theta_{\tilde m,\mu})_{br}$ by uniqueness of the steady-state (see Section \ref{sec:motiv}). Hence, $\theta_{\tilde m,\mu}$ is non-decreasing. Moreover, according to \eqref{LDE}, $\theta_{\tilde m,\mu}$ is convex on $(0,1-\ell)$ and concave on $(1-\ell,1)$ which, combined with the boundary conditions on $\theta_{\tilde m,\mu}$, justifies the positiveness of its derivative. The expected result follows.

%%%Convergence des suites de maximiseurs.
\paragraph{Step 1: convergence of sequences of maximizers.}
As a consequence of Theorem \ref{Concentration}, we get  that the functions $\tilde m=\kappa\chi_{(0,\ell)}$ or $\tilde m(1-\cdot)=\kappa \chi_{(1-\ell,1)}$ are the only closure points of the family $(m_\mu)_{\mu>0}$ for the $L^1(0,1)$ topology. 
\paragraph{Step 2: asymptotic behaviour of $p_{m,\mu}$ and of $\p_{m,\mu}$}
We claim that, as done for the solution $\theta_{m,\mu}$ of \eqref{LDE}, the following asymptotic behaviour for the adjoint state $p_{m,\mu}$
$$
p_{m,\mu}=-\frac1{m_0}+\underset{\mu\to\infty}{\operatorname{O}}(\frac1\mu),\quad \text{in }W^{2,2}(0,1),
$$
by using Sobolev embeddings. In particular, this expansion holds in $\mathscr C^1([0,1])$.
\\ Introduce the function $z_\mu=\mu(\varphi_{m_\mu,\mu}+1)$. Using the convergence results established in the previous steps, in particular that $(m_\mu)_{\mu>0}$ converges to $\tilde m$ in $L^1(0,1)$ and that $\varphi_{m_\mu,\mu}=-1+\underset{\mu\to\infty}{\operatorname{O}}(\frac1\mu)$ uniformly in $\mathscr C^{1,\alpha}([0,1])$\footnote{This is obtained similarly to the proof's technique of theorem \ref{TheoremePrincipal}, using elliptic estimates and Sobolev embedding for the functions $\theta_{m,\mu}$ and $p_{m,\mu}$. } as $\mu\to +\infty$, one infers that $(z_\mu)_{\mu>0}$
 is uniformly bounded in $\mathscr C^{1,\alpha}([0,1])$ and converges, up to a subsequence to $z_\infty$ in $\mathscr C^1([0,1])$, where $z_\infty$ satisfies in particular
 $$
 z_\infty''+2(m_0-\tilde m)=0,
 $$
 with Neumann Boundary conditions  in the $W^{1,2}$ sense. 
\paragraph{Conclusion: $m_\mu=\tilde m$ or $\tilde m(1-\cdot)$ whenever $\mu$ is large enough.}
According to Theorem \ref{TheoremePrincipal} and Proposition \ref{propoptimCions}, we know at this step that for $\mu$ large enough, there exists $c_\mu\in \R$ such that 
 $$
 \{\varphi_{m_\mu,\mu}> c_\mu\}=\{m_\mu=0\},\qquad \{\varphi_{m_\mu,\mu}<c_\mu\}=\{m_\mu=\kappa\}.
 $$

 We will show that, provided that $\mu$ be large enough, one has necessarily $m_\mu=\tilde m$ or $m_\mu=\tilde m(1-\cdot)$.
 Since $\tilde m= \kappa$ in $(0,\ell)$, it follows that $z_{\infty}$ is strictly convex on this interval and since $z_{\infty}'(0)=0$, one has necessarily $z_{\infty}'>0$ in $(0,\ell)$. Similarly, by concavity of $z_\infty$ in  $(\ell,1)$, one has 
 $z_{\infty}'>0$ in this interval.
 
Furthermore, let us introduce $d_\mu=\mu(c_\mu+1)$. Since $(z_\mu)_{\mu>0}$ is bounded in $C^0((0,1))$, $(d_\mu)_{\mu>0}$ converges up to a subsequence to some $d_\infty$. By monotonicity of $z_\infty$ and a compactness argument, there exists a unique $x_{\infty}\in [0,\ell]$ such that $z_{\infty}(x_{\infty})=d_{\infty}$. The dominated convergence theorem hence yields 
$$
| \{z_\infty\leq  d_\infty\}| = \kappa \ell, \qquad | \{z_\infty\geq  d_\infty\} |= \kappa (1-\ell),
$$
and the the aforementioned local convergence results yield 
$$
\{z_\infty>d_\infty\}\subset\{\tilde m=0\}\, , \quad \{z_\infty<d_\infty\}\subset \{\tilde m=\kappa\}.
$$
Hence, the inclusions are equalities (the equality of sets must be understood up to a zero Lebesgue measure set) by using that $z_{\infty}$ is increasing. 

Moreover, since $z_{\infty}$ is increasing, one has $z_{\infty}(0)<d_{\infty}$ and $z_{\infty}(1)>d_{\infty}$. 
Since the family $(z_{\mu})_{\mu>0}$ is uniformly Lipschitz-continuous, there exists $\e>0$ such that for $\mu$ large enough, there holds
$$
z_{\mu} <d_{\mu} \hbox{ in } (0,\e), \quad z_{\mu} >d_{\mu} \hbox{ in } (1-\e,1), \quad 
z_{\mu}'>0 \quad \hbox{ in } (\e, 1-\e).
$$

This implies the existence of $x_{\mu}\in (0,1)$ such that 
$$
\{z_{\mu}<d_{\mu}\}=[0,x_{\mu}) \quad \hbox{and} \quad \{z_{\mu}>d_{\mu}\}=(x_{\mu},1],
$$
whence the result.

\subsection{Proof of Theorem \ref{fragmentation}}
\color{black}Let $\kappa>0, m_0>0$, and $\tilde m:=\kappa \chi_{[1-\ell,1)}$ with $\ell=\frac{m_0}\kappa$, i.e the single crenel distribution.\color{black}
\\In order to prove this result, as the function $\mu>0\mapsto {\mathcal F_\mu} \left(\tilde m(2~\cdot)\right)$   
has a  first   local maximizer  (\cite[Theorem 1.2, Remark 1.4]{LouMigration}), we  define $\mu_{1}$ as its first local maximizer. One gets from a simple change of variables that 
$\theta_{\tilde m,\mu_{1}}(2 x )=\theta_{\tilde m(2\cdot ),\mu_{1}/4}(x)$ for all $x\in \O$ and thus one has 
$$
F_{\mu_1}(\tilde m)={F}_{\frac{\mu_1}4}\left(\tilde m(2\cdot)\right)
$$
But our choice of $\mu_{1}$  yields that $\mu\mapsto {\mathcal F_\mu} \left(\tilde m(2\cdot)\right) $ is increasing on $(0,\mu_{1})$ and thus:
\begin{equation}
F_{\mu_1}(\tilde m)={F}_{\frac{\mu_1}4}\left(\tilde m(2\cdot)\right)<{F}_{\mu_1} \left(\tilde m(2\cdot)\right).
\end{equation}

\section{Conclusion and further comments}
\subsection{About the 1D case} \label{sec:1D}
\begin{figure}[h]
\begin{center}
\includegraphics[width=5.5cm]{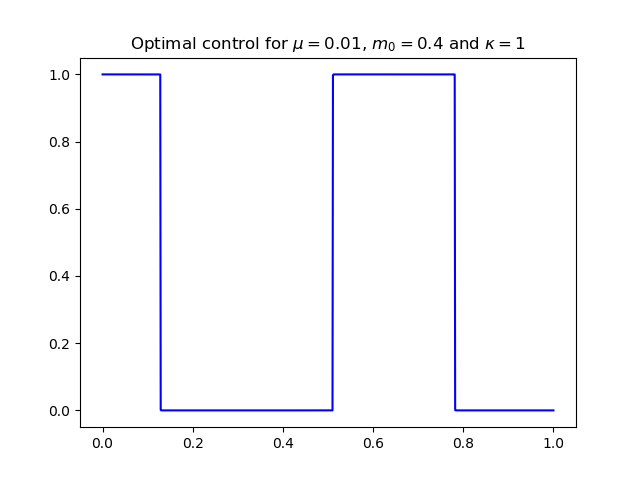}\hspace{-1cm}
\includegraphics[width=5.5cm]{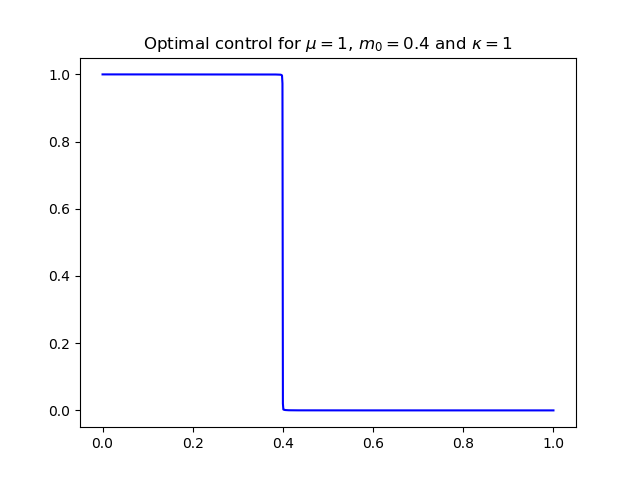}\hspace{-1cm}
\includegraphics[width=5.5cm]{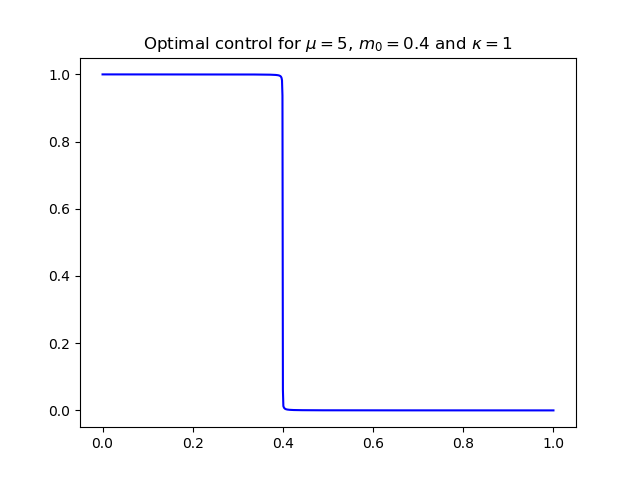}\\
\includegraphics[width=5.4cm]{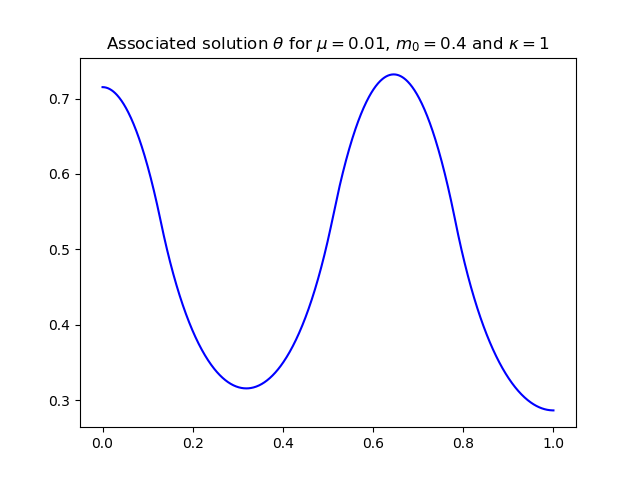}\hspace{-0.7cm}
\includegraphics[width=5.4cm]{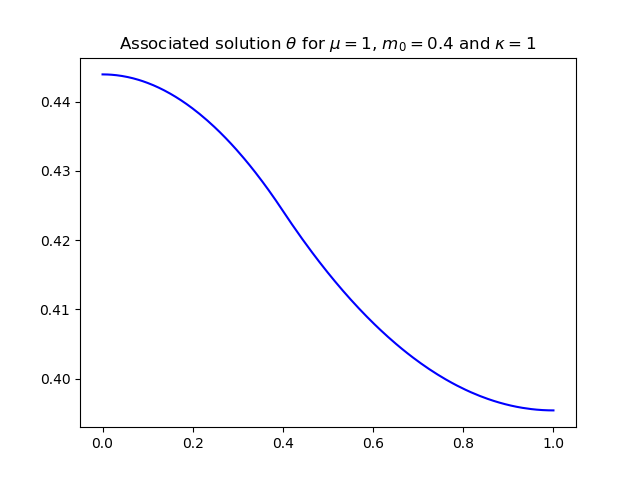}\hspace{-1cm}
\includegraphics[width=5.4cm]{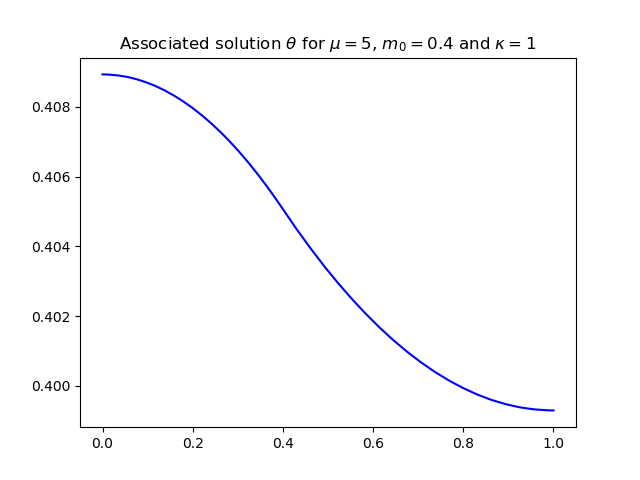}
\caption{$m_0=0.4$, $\kappa=1$. From left to right: $\mu=0.01,1,5$. Top: plot of the optimal solution of Problem \eqref{TSOM} computed with the help of an interior point method. Bottom: plot of the corresponding eigenfunction.  \label{figint2}}
\end{center}
\end{figure}

Let us assume in this section that $n=1$ and $\Omega=(0,1)$. We provide hereafter several numerical simulations based on the primal formulation of the optimal design problem \eqref{TSOM}: 
%on Fig. \ref{figint1}, the objective function is plotted with respect to $x_0$ for several values of $\mu$, where we assumed the control function $m$ having the particular form $m=\kappa \chi_I$ with $I=(x_0-m_0/2\kappa,x_0+m_0/2\kappa)$. Here, $x_0$ is the middle of the interval. 
on  Fig. \ref{figint2}, we investigate the general problem \eqref{TSOM} and we plot the optimal $m$ determined numerically for several values of $\mu$.

These simulations were obtained with an interior point method applied to the optimal control problem \eqref{TSOM}. We used a Runge-Kutta method of order 4 to discretize the underlying differential equations. The control $m$ has been also discretized, which has allowed to reduce the optimal control problem to some finite dimensional minimization problem with constraints. We used the code \texttt{IPOPT} (see \cite{IPOPT}) combined with \texttt{AMPL} (see \cite{AMPL}) on a standard desktop machine. We considered a regular subdivision of $(0,1)$ with $N$ points, where the order of magnitude of $N$ is $1000/\mu$. The resulting code works out the solution quickly (around 5 to 10 seconds depending on the choice of the parameter $\mu$).

%\begin{figure}[h]
%\begin{center}
%\includegraphics[width=5.5cm]{int1m040kappa1mu05.eps}\hspace{-1cm}
%\includegraphics[width=5.5cm]{int1m040kappa1mu1.eps}\hspace{-1cm}
%\includegraphics[width=5.5cm]{int1m040kappa1mu5.eps}
%\caption{$m_0=0.4$, $\kappa=1$. Graph of $\int_0^1\theta_{m,\mu}$ with respect to $x_0$ where $m=\kappa\chi_{(x_0-m_0/2,x_0+m_0/2)}$. From left to right: $\mu=0.5,1,5$.\label{figint1}}
%\end{center}
%\end{figure}
%

%These simulations highlighted that, in the 1D case, assuming that $m=\kappa \chi_I$, where $I$ is an interval, the best locations of $I$ are the extreme one (so that $m$ is either decreasing or increasing on $\O$).
\red{
In the cases mentioned above, the algorithm is initialized with several choices of function $m$, among which 
the optimal simple crenel as $\mu$ is large enough. If $\mu$ is equal to 1 or 5, the simple crenel is obtained at convergence. Nevertheless, in the case $\mu=0.01$, we obtain a ``symmetric'' double crenel (in accordance with Theorem \ref{fragmentation}) at convergence.}

\red{Although we have no guarantee to obtain optimal solutions by using this numerical approach, we checked that a simple crenel is better than a double one in the cases $\mu=1, \ 5$ whereas we observe the contrary in the case  $\mu=0.01$. }
 
\red{Notice that we encountered a problem when dealing with too small values of $\mu$ (for instance $\mu=0.001$). Indeed, in that case, the stiffness of the discretized system seems to become huge as $\mu$ takes small positive values and makes the numerical computations hard to converge. Improvements of the numerical method should be found for further numerical investigations. }

%\subsection{Comments on the Dirichlet case}
%When dealing with Dirichlet boundary conditions instead of Neumann ones, all the results of this article can be directly extended to that case, by noting that the assumption $\lambda_1(m,\mu)>0$ needs to be replaced by the condition
%$$
%\lambda_1(m,\mu)>\lambda_1(0,\mu)+m_0,
%$$
%that is,  we need to replace \eqref{H1} with \eqref{H1'}:
%\begin{equation}\label{H1'}\tag{H1'}
%m_0+\lambda_1(0,\mu)>0.
%\end{equation}
%Once this condition is satisfied, all of the previous proofs can be written along the same lines and are even sometimes simplified. Indeed, one of the main difficulties when dealing with Neumann boundary conditions  was to control the integral quantities $\eta_{k,m}$ by using a Poincar\'e-Wirtinger inequality. Working with Dirichlet boundary conditions enables to work with the standard Poincar\'e inequality in $W^{1,2}_0$ and to overcome this technical point.

%It is also interesting to note that, in dimension $N\geq 2$ with Dirichlet boundary conditions and $\O=\mathbb{B}(0;R)$, the approach developed within this article can be adapted to prove the existence of $\mu^*$ such that, for any $\mu\geq \mu^*$, any sequence of maximizers of ${\mathcal F_\mu}$ converges in $L^1(\O)$ to $m^*$, where $m^*=\kappa \chi_{\mathbb{B}(0,r)}$ and $r$ is chosen so that $\int_\O m^*=m_0|\O|$.
%Note that the symmetrizations used in the proof of Theorem \ref{Creneau} have to be replaced by the so-called (radial) Schwarz symmetrization (see e.g. \cite{MR810619,henrot-pierre}).

\subsection{Comments and open issues}
It is also interesting, from a biological point of view, to investigate a more general version of Problem \eqref{TSOM} for changing-sign weights. In that case, the admissible class of weights is then transformed (for instance) into 
$$
\widetilde{\mathcal{M}}_{m_0,\kappa}(\O)=\left\{ m\in L^\infty(\O)\, ,m\in [-1;\kappa]\text{ a.e and }\fint_\O m=m_0\right\},
$$
with $m_0\in (0,1)$ (so that $\lambda_1(m,\mu)>0$ and Equation \eqref{LDE} is well-posed). We claim that the main results of this article can be extended without effort to this new framework and that we will still obtain the bang-bang character of maximizers provided that $\mu$ be large enough.
Such a class has also been considered in the context of principal eigenvalue minimization (see \cite{HinKaoLau12,LamboleyLaurainNadinPrivatProperties}).

Finally, we end this section by providing some open problems for which we did not manage to bring complete answer and that deserve and remain, to our opinion, to be investigated. They are in order:
\begin{itemize}
\item (for general domains $\Omega$) we conjecture that maximizers are bang-bang functions for any $\mu>0$. As outlined in the introduction, this conjecture is supported by Theorem \ref{TheoremePrincipal} and the main result of \cite{NagaharaYanagida}.
\item (for general domains $\Omega$) use the main results of the present article to determine numerically the maximizer $m^*$ with the help of an adapted shape optimization algorithm;  
\item (for $\O=(0;1)$) given that, for $\mu$ small enough, the optimal configurations for $\lambda_1(\cdot,\mu)$ and ${\mathcal F_\mu}$ are not equal, it would be natural  and biologically relevant to try to maximize a convex combination of ${\mathcal F_\mu}$ and $\lambda_1(\cdot,\mu)$.
\item  (for general domains $\Omega$) investigate the asymptotic behavior of maximizer as the parameter $\mu$ tends to 0? Such a issue appears intricate since it requires a refine study of singular limits for Problem \eqref{LDE}.
\end{itemize}
\section*{Aknowledgements}
The authors thank the reviewers for their thorough reviews and highly appreciate the comments and suggestions, which significantly contributed to improve the quality of the publication.

\appendix
\section{Convergence of the  series}\label{append:CVSeries}

Let $\frac{1}{\mu^*_1}$ be the minimum of the convergence radii associated to the power series $\sum \alpha(k)x^k$, $\sum \sigma(k)x^k$, $\sum \gamma(k)x^k$, $\sum \delta(k)x^k$ and $\sum \varepsilon(k)x^k$ introduced in the proof of Theorem \ref{TheoremePrincipal}.

We will show that, whenever $\mu\geq\mu^*_1$, the following expansions
$$
\sum_{\ell=0}^{+\infty}\frac{\etakm}{\mu^k}=\tmm, \quad \sum_{k=1}^{+\infty}\frac{\detakm}{\mu^k}=\dtm,\quad  \sum_{k=1}^{+\infty}\frac{\ddot{\eta}_{k,m}}{\mu^k}=\ddot{\theta}_{m,\mu}
$$
make sense in $\l$. Since the proofs for the series defining $\dtm$ and $\ddot\theta_{m,\mu}$ are exactly similar to the one for  $\tmm$, we only concentrate on the expansion of $\tmm$. By construction, the series $g_{\infty,\mu}:=\sum_{\ell=0}^{+\infty}\frac{\etakm}{\mu^k}$ converges in $W^{1,2}(\O)$ to a function $g_{\infty,\mu}$. We need to show that $g_{\infty,\mu}=\tmm$.

To this aim, let us set
$$
g_{N,\mu}:=\sum_{k=0}^N \frac{\etakm}{\mu^k}
$$
 for any $N\in \N^*$, 
Notice that $g_{N,\mu}$ solves the equation
\begin{equation}\label{Eq:gN} \mu \Delta g_{N,\mu}+g_{N-1,\mu}m-\sum_{k=0}^N \frac{\eta_{k,m}}{\mu^k}g_{N-k,\mu}=0,\quad \text{in }\O
\end{equation}
with Neumann boundary conditions.

In order to pass to the limit $N\to \infty$, one has to determine the limit of $\tilde g_{N,\mu}:=\sum_{k=0}^N \frac{\eta_{k,m}}{\mu^k}g_{N-k,\mu}$.  First note that the Cauchy-Schwarz inequality proves the absolute convergence of the sequence $\left\{\tilde g_{N,\mu}\right\}_{N\in \N}$ in $W^{1,2}(\O)$ as $N\to \infty$. Let $H$ denote its limit.
Now, let us show that 
$$
\tilde g_{N,\mu}\underset{N\to \infty}\to g_{\infty,\mu}^2\quad \text{ in }L^2(\O),
$$  
whenever $\mu$ is large enough. Let $R_1$ be the convergence radius of the power series associated with the sequence $\{\alpha(k)\}_{k\in \N}$. This convergence radius is known to be positive. As a consequence, the convergence radius $R_2$ of the power series associated with the sequence $\{\alpha(k)^2\}_{k\in \N}$ is also positive and $R_2=R_1^2.$

Let $\e>0$. Since we are only working with large diffusivities, let us assume that $\mu\geq 1$ and  that $\mu>(\frac1{R_2})^{1/\e}.$
Noting that, for any $N\in \N$, we have
\begin{align*}
\tilde g_{N,\mu}-g_{N,\mu}^2&=\sum_{k=0}^N \frac1{\mu^k}\eta_{k,m}\left( \sum_{\ell=N-k+1}^N\frac{\eta_{\ell,m}}{\mu^\ell}\right).
\end{align*}
and using the fact that the sequence $\{\alpha(k)\}_{k\in \N}$ was built increasing, we get the existence of $M>0$ such that 
\begin{align*}
\frac{1}{|\Omega|}\Vert \tilde g_{N,\mu}-g_{N,\mu}^2\Vert_{L^2(\O)}&=\frac{1}{|\Omega|}\Bigg\Vert\sum_{k=0}^N \frac1{\mu^k}\eta_{k,m}\left( \sum_{\ell=N-k+1}^N\frac{\eta_{\ell,m}}{\mu^\ell}\right)\Bigg\Vert_{L^2(\O)}&
\\&\leq \sum_{k=0}^N\frac{\alpha(k)}{\mu^k}\left(\sum_{\ell=N-k+1}^N\frac{\alpha({\ell})}{\mu^\ell}\right)&
\\&\leq \alpha(N)^2\sum_{k=0}^N\frac1{\mu^k}\left(\frac1{\mu^{N-k+1}}\frac{1-\frac{1}{\mu^{k-1}}}{1-\frac1\mu}\right) &\text{ by using that } \alpha(k)\alpha(\ell)\leq \alpha(N)^2
\\&\leq M\frac{(N+1)\alpha(N)^2}{\mu^{N+1}}
\\&=M\frac{N+1}{(\mu^{1-\e})^{N+1}}\frac{\alpha(N)^2}{(\mu^\e)^{N+1}}.
\end{align*}
This last quantity converges to zero as $N\to \infty$. Besides, since $\mu^\e\geq \frac1{R_2}$, it follows that the sequence $\left\{\frac{\alpha(N)^2}{(\mu^\e)^{N+1}}\right\}_{N\to \infty}$ is bounded. Assuming moreover that $\mu^{1-\e}>1$, we get
$$
\frac{N+1}{(\mu^{1-\e})^{N+1}}\xrightarrow[N\to \infty]{} 0.
$$ 
We conclude that $H=g_{\infty,\mu}^2$. 
Passing to the limit in Equation \eqref{Eq:gN}, it follows that 
$$
\mu \Delta g_{\infty,\mu}^2+g_{\infty,\mu}(m-g_{\infty,\mu})=0\quad \text{in }\O
$$
 with Neumann boundary conditions.\\  
 
 Finally, we know that $g_{\infty,\mu}\underset{\mu \rightarrow +\infty}\rightarrow m_0$ uniformly in $\mathcal M_{m_0,\kappa}(\O)$ and moreover, one has $m_0>0$. It follows that, for $\mu$ large enough, $g_{\infty,\mu}$ is positive. The uniqueness of positive solutions of equation \eqref{LDE} entails that, for $\mu$ large enough, $g_{\infty,\mu}=\tmm$.
This concludes the proof of the series expansion convergences and thus, the proof of Theorem \ref{TheoremePrincipal}.

\bibliographystyle{abbrv}
\nocite{*}
\bibliography{biblio}
\addcontentsline{toc}{part}{Bibliography}

\end{document}